\documentclass[reqno]{nyjm}

\usepackage{general,resistance,cancel,ifthen ,leading}
\usepackage[switch*,pagewise]{lineno}

\usepackage[bookmarks, colorlinks=true, pdfstartview=FitV, linkcolor=blue, citecolor=blue, urlcolor=blue]{hyperref}
\usepackage{cite}  % <= put after hyperref

\usepackage[left=1.85in,right=1.85in,top=1.7in,bottom=1.40in]{geometry}%showframe

\numberwithin{equation}{section} \numberwithin{theorem}{section}

\renewcommand{\cpath}{\gw}     %path (sample) of the random walk
\renewcommand{\linenopax}{}

\begin{document}
  %\linenumbers
  %\pagewiselinenumbers

%\begin{frontmatter}
\title[Gel'fand triples and boundaries of infinite networks]{Gel'fand triples and boundaries \\ of infinite networks}

%    author one information
\author{Palle E. T. Jorgensen}
\address{University of Iowa, Iowa City, IA 52246-1419 USA}
\email{jorgen@math.uiowa.edu}

%    author two information
\author{Erin P. J. Pearse}
\address{University of Oklahoma, Norman OK 73019-0315 USA}
\email{ep@ou.edu}

\thanks{The work of PETJ was partially supported by NSF grant DMS-0457581. The work of EPJP was partially supported by the University of Iowa Department of Mathematics NSF VIGRE grant DMS-0602242.}

\begin{abstract}
  We study the boundary theory of a connected weighted graph $G$ from the viewpoint of stochastic integration. For the Hilbert space \HE of Dirichlet-finite functions on $G$, we construct a Gel'fand triple $S \ci {\mathcal H}_{\mathcal E} \ci S'$. This yields a probability measure $\mathbb{P}$ on $S'$ and an isometric embedding of ${\mathcal H}_{\mathcal E}$ into $L^2(S',\mathbb{P})$, and hence gives a concrete representation of the boundary as a certain class of ``distributions'' in $S'$. In a previous paper, we proved a discrete Gauss-Green identity for infinite networks which produces a boundary representation for harmonic functions of finite energy, given as a certain limit.
In this paper, we use techniques from stochastic integration to make the boundary $\operatorname{bd}G$ precise as a measure space, and obtain a boundary integral representation as an integral over $S'$.  
\end{abstract}
 %prove a discrete Gauss-Green identity constructed a reproducing kernel $\{v_x\}$ 
 %\[{\mathcal E}(u,v) = \sum_{G} u \Delta v + \sum_{\operatorname{bd}G} u \frac{\partial}{\partial \mathbf{n}} v,\]
 %where the latter sum is understood in a limiting sense. Applying this formula to a harmonic function $u \in {\mathcal H}_{\mathcal E}$ gives a boundary representation



%  \dedicatory{Dedicated to Bob Powers, in recognition of his pioneering spirit.}

  \date{\today.}

\maketitle

\setcounter{tocdepth}{1}
%\version{}{\setcounter{tocdepth}{2}}
{\small \tableofcontents}

  \leading{14pt} %12pt
\allowdisplaybreaks

%%!TEX root = bdG.tex

\section{Introduction}
\label{sec:introduction}

In this paper, we develop a boundary theory for an infinite network (connected weighted graph) \Graph, using some techniques from the theory of stochastic integration. For the Hilbert space \HE of finite-energy functions on $G$, we construct a Gel'fand triple $\Schw \ci \HE \ci \Schw'$, where both containments are strict, and the inclusion mappings are continuous. Here, \Schw is a space of ``test functions'' on the network and $\Schw'$ is a class of ``distributions'' on the network, analogous to Schwartz's classical \emph{functions of rapid decay} and \emph{tempered distributions}, respectively.
A key result of this paper is Theorem~\ref{thm:HE-isom-to-L2(S',P)}, which establishes an isometric embedding of \HE into the Hilbert space $L^2(\Schw',\prob)$, where \prob is a Gaussian probability measure on $\Schw'$.
In Martin boundary theory, elements of the boundary are understood as certain minimal harmonic functions. We are studying finite-energy harmonic functions instead of positive harmonic functions, but the construction outlined above allows us to study elements of the boundary in an analogous fashion.

In a previous paper, we proved a discrete Gauss-Green identity for infinite networks:
  \linenopax
  \begin{equation}\label{eqn:DGG-preview}
    \energy(u,v)
    = \sum_{\verts} \cj{u} \Lap v
      + \sum_{\bd \Graph} \cj{u} \dn v.
  \end{equation} 
Here, \energy is a Dirichlet form and \Lap is the graph Laplacian; see Theorem~\ref{thm:Discrete-Gauss-Green-identity} and the discussion preceding it for precise definitions of the other terms.
  Formula \eqref{eqn:DGG-preview} yields a boundary representation for a harmonic function $u$: 
  \linenopax
  \begin{align}\label{eqn:bdy-repn-preview}
    u(x) = \sum_{\operatorname{bd}G} u \frac{\partial v_x}{\partial \mathbf{n}} + C, 
  \end{align}
where $C$ is a constant and the sum is actually defined as a limit of ``Riemann sums'' over an increasing sequence of finite subnetworks of \Graph; see Definition~\ref{def:boundary-sum}.
In this paper, we use functional integration techniques from stochastic integration to make the boundary $\operatorname{bd}G$ precise as a measure space, and replace the sum with an integral over $\Schw'$, thus obtaining a boundary integral representation (in a sense analogous to that of Poisson or Martin boundary theory) for the harmonic function $u$.  

%As a corollary (Corollary~\ref{thm:Boundary-integral-repn-for-harm}), we obtain a (relatively concrete) boundary representation for the harmonic functions of finite energy. 
Another key result of this paper is Corollary~\ref{thm:Boundary-integral-repn-for-harm}, which follows readily from Theorem~\ref{thm:HE-isom-to-L2(S',P)} and gives a boundary integral representation for harmonic functions of finite energy:
  \begin{equation}\label{eqn:integral-boundary-repn-of-h-preview}
    u(x) = \int_{\Squoth} u(\gx) h_x(\gx) \, d\prob(\gx) + u(o).
  \end{equation} 
Here, $\{h_x\}_{x \in \Graph}$ is a family of harmonic functions parametrized by the vertices and discussed in detail just below (and see also Lemma~\ref{thm:MP-contains-spans}). Given a transient network, this allows one to identify a space of functions in $\Schw'$ corresponding to the boundary of the network (in a manner reminiscent of the Martin boundary). 
Additionally, Example~\ref{exm:a,b-ladder} presents the construction of a harmonic function of finite energy on a network with one ``graph end'' (in fact, a two-parameter family of such networks). The existence of such functions was first proved in \cite{CaW92}, but we have never seen an explicit formula given before. We now proceed to describe these results in a bit more detail. The reader is also referred to the survey paper \cite{RBIN} which gives an overview of how the results of the present paper fit into a larger investigation of functions of finite energy on resistance networks, and the effective resistance metric.

\subsection{Overview}

We define what it means for a function on a network to have finite energy in Definition~\ref{def:graph-energy}. Then the (discrete) Dirichlet energy form $\mathcal E$ (also given in Definition~\ref{def:graph-energy}) is an inner product on the space of functions of finite energy, and in fact produces a Hilbert space which we denote by \HE. The space \HE consists of equivalence classes of functions on the vertices of \Graph, where $u \simeq v$ iff $u-v$ is a constant function. In a previous paper, we constructed a reproducing kernel $\{v_x\}_{x \in \verts}$ for this Hilbert space and used it to prove a discrete Gauss-Green identity which is recalled in Theorem~\ref{thm:Discrete-Gauss-Green-identity}. 

The space \HE also enjoys an orthogonal decomposition into the subspace \Fin of (\energy-limits of) finitely supported functions and the subspace \Harm of harmonic functions; see Definitions~\ref{def:Fin}--\ref{def:Harm} and Theorem~\ref{thm:HE=Fin+Harm}. Since the reproducing kernel behaves well with respect to projections, we also have reproducing kernels $\{f_x\}_{x \in \verts}$ for \Fin and $\{h_x\}_{x \in \verts}$ for \Harm, where
\linenopax
\begin{align*}%\label{eqn:repkernels-for-Fin-and-Harm}
  f_x := \Pfin v_x,
  \q\text{and}\q
  h_x := \Phar v_x.
\end{align*}
It should be noted that these kernels are reproducing up to an additive constant; in other words, for some fixed reference vertex $o$,
\linenopax
\begin{align*}%\label{eqn:}
  \la v_x, u \ra_\energy = u(x) - u(o),
\end{align*}
and similarly for $f_x$ and $h_x$.

Recall the classical result of Poisson that gives a kernel $k:\gW \times \del \gW \to \bR$ from which a bounded harmonic function can be given via
\linenopax
  \begin{equation}\label{eqn:Poisson-bdy-repn}
    u(x) = \int_{\del \gW} u(y) k(x,dy),
    \qq y \in \del \gW.
  \end{equation}
We are motivated by the discrete analogue of this result appearing in \eqref{eqn:bdy-repn-preview}.

%\begin{theorem}[{\cite[Cor.~3.14]{DGG}}]
%  \label{thm:boundary-repn-for-harmonic}
%  For all $u \in \Harm$, and $h_x = \Phar v_x$,
%  \linenopax
%  \begin{equation}\label{eqn:boundary-repn-for-harmonic}
%    u(x) - u(o) = \sum_{\bd \Graph} u \dn{h_x}.
%  \end{equation}
  %\begin{proof}
  %  Compute $u(x)-u(o) = \la h_x, u \ra_\energy = \cj{\la u, h_x \ra_\energy} = \sum_{\bd \Graph} u \dn{h_x}$ because $\sum_{\verts} u \Lap h_x = 0$. (By Lemma~2.24, $h_x$ is \bR-valued.)
  %\end{proof}
%\end{theorem}

%This gives a boundary sum representation of harmonic functions, but the boundary sum in \eqref{eqn:boundary-repn-for-harmonic} is understood only as a limit of sums. 
Comparison of \eqref{eqn:bdy-repn-preview} and \eqref {eqn:Poisson-bdy-repn} makes one optimistic that $\bd \Graph$ can be realized as a measure space which supports a measure corresponding to $\dn{h_x}$, thus replacing the sum in \eqref{eqn:bdy-repn-preview} with a integral. In Corollary~\ref{thm:Boundary-integral-repn-for-harm}, we do precisely this.
  
 Boundary theory of harmonic functions can roughly be divided three ways: the bounded harmonic functions (Poisson theory), the nonnegative harmonic functions (Martin theory), and the finite-energy harmonic functions studied in the present paper. While Poisson theory is a subset of Martin theory, the relationship between Martin theory and the study of \HE is more subtle. For example, there exist unbounded functions of finite energy; cf.~\cite[Ex.~13.10]{OTERN}. Some results detailing the interrelations  are given in \cite[\S3.7]{Soardi94}.
  
  Whether the focus is on the harmonic functions which are bounded, nonnegative, or finite-energy, the goals of the associated boundary theory are essentially the same:
\begin{enumerate}
  \item Construct a space $\cj{\sD}$ which extends the original domain \sD; this can be done by taking closure, compactifying, or similar operations.
  \item One can then identify the boundary $\bd \sD$ as $\cj{\sD} \less \sD$, or (if the boundary thus obtained would be larger than necessary/practical for the application in mind), as some subset of $\cj{\sD} \less \sD$.
  \item Define a procedure for extending harmonic functions $u$ from \sD to $\bd \sD$. This extension $\tilde u$ may be a measure (or other linear functional) on $\bd \sD$; it may not be a function.
  \item Obtain a kernel $\mathbbm{k}(x,\gb)$ defined on $\sD \times \bd \sD$ against which one can integrate the extension $\tilde u$ so as to recover the value of $u$ at a point in \sD:
  \linenopax
  \begin{align*}%\label{eqn:}
    u(x) = \int_{\bd \sD} \mathbbm{k}(x,\gb) \tilde u(d\gb),
    \q \forall x \in \sD,
  \end{align*}
  whenever $u$ is a harmonic functions of the given class.
\end{enumerate}
Our approach to (1) is to use Gel'fand triples to extend the original domain, a method which is novel as far as we know. In a forthcoming work \cite{Shannon}, we will introduce an interpolation formula that uses the analytic framework developed in this paper, and which turns \Graph into a stochastic process. In particular, the interpolation formula allows one to find continua which naturally extend \Graph. 
  
The difference between our boundary theory and that of Poisson and Martin is rooted in our focus on \HE rather than $\ell^2$: both classical theories concern harmonic functions with growth/decay restrictions. By contrast, provided they neither grow too wildly nor oscillate too wildly, elements of \HE may have values tending to both $+\iy$ and $-\iy$. See \cite[Ex.~13.10]{OTERN} for a function $h \in \Harm$ which is unbounded in this way. \emph{Positive} harmonic functions are naturally given to analysis based on probabilistic and potential-theoretic techniques, and the companion study of superharmonic (or subharmonic) functions is indispensable. Without positivity, however, one can get more mileage by considering the Dirichlet form \energy as an inner product and studying the resulting Hilbert space geometry.

Our boundary essentially consists of (equivalence classes of) infinite paths which can be distinguished by monopoles, i.e., two paths are not equivalent iff there is a monopole $w$ with different limiting values along each path. It is an immediate consequence that recurrent networks have no boundary, and transient networks with no nontrivial harmonic functions have exactly one boundary point (corresponding to the fact that the monopole at $x$ is unique). In particular, the integer lattices $(\bZd,\one)$ each have 1 boundary point for $d \geq 3$ and 0 boundary points for $d=1,2$. In contrast, the Martin boundary of $(\bZ^d,\one)$ is homeomorphic to the unit sphere $S^{n-1}$ (where $S^0 = \{-1,1\}$), and each $(\bZ^d,\one)$ has only one graph ends (except for $(\bZ,\one)$, which has two); cf.~\cite[\S3.B]{PicWoess90}, for example.

\subsection{Outline}

In our version of the program outlined above, we follow the steps in the order (2)-(3)-(1). A brief summary is given here; further introductory material and technical details appear at the beginning of each subsection. 

\S\ref{sec:electrical-resistance-networks} recalls basic definitions and some previously obtained results.

\S\ref{sec:Gel'fand-triple-for-HE} describes two methods for constructing a Gel'fand triple.
The technique presented in \S\ref{sec:Gel'fand-triples-via-Gram-Schmidt} works for any network $(\Graph,\cond)$ and makes use of an orthonormal basis of \HE derived from the energy kernel $\{v_x\}_{x \in \Graph}$ via the Gram-Schmidt algorithm, or equivalently, from the domain of a certain operator \Onb.
The approach given in \S\ref{sec:unbounded-case} works only for networks where \Lap is an unbounded operator on \HE. This version of \Schw is constructed in terms of the domain of \Lap. %, and Ito integration is used to extend the resultant pairing on $\Schw \times \Schw'$ to all of $\HE \times \Schw'$. We also obtain some more specific information about the structure of $\Schw'$ in this case. %This yields a suitable class of linear functionals \gx on \HE, and we can extend a function $u$ on \HE to $\tilde u$ on $\Schw'$ by duality, i.e., $\tilde u(\gx) := \la u, \gx\ra_\Wiener$. 
%We need to expand the scope of enquiry to include $\Schw'$ because \HE will not be sufficient; no infinite-dimensional Hilbert space can support a \gs-finite probability measure, by a theorem of Nelson.

\S\ref{sec:Structure-of-Schw} studies the structure of \Schw (the space of test functions) and $\Schw'$ (the space of distributions) and establishes some key lemmas for later use.

\S\ref{sec:Wiener-embedding} proves a key result: Theorem~\ref{thm:HE-isom-to-L2(S',P)}, which establishes %an isometric embedding of the space of functions of finite energy (defined on vertices of a network) into a Hilbert space $L^2(\Schw',\prob)$ 
the isometric embedding of \HE into $L^2(\Schw',\prob)$ given by the Wiener transform. Applying this isometry to the energy kernel $\{v_x\}$, we get a reproducing kernel $\mathbbm{k}(x,d \prob)$ given in terms of a version of Wiener measure. In fact, \prob is a Gaussian probability measure on $\Schw'$ whose support is disjoint from \Fin. The results in this section hold for any Gel'fand triple; in particular, for either of the ones constructed in \S\ref{sec:Gel'fand-triple-for-HE} and \S\ref{sec:unbounded-case}.

\S\ref{sec:Operator-theoretic-interpretation-of-bdG}
 We consider certain measures $\gm_x$, defined in terms of the kernel and the Wiener measure just introduced, which are supported on $\Schw'/\Fin$ and indexed by the vertices $x \in \verts$. Then points of $\bd G$ correspond to limits of sequences $\{\gm_{x_n}\}$ where $x_n \to \iy$, modulo a suitable equivalence relation. %This is the content of \S\ref{sec:bdG-as-equivalence-classes-of-paths}.

\version{}{\S\ref{sec:examples} presents some elementary but illuminating examples.}

\begin{remark}\label{rem:comparison-to-path-space-measures}
  %In \S\ref{sec:Probabilistic-interpretation} we will return to the three-way comparison of harmonic functions which are bounded, nonnegative, or finite-energy, but for a different purpose: the construction of measures on spaces of (infinite) paths in $(\Graph, \cond)$. In the case of bounded harmonic functions on $(\Graph, \cond)$, the associated probability space is derived directly as a space of infinite paths in \Graph, and the measure is constructed via the standard Kolmogorov consistency method. That is, as a projective limit constructed via cylinder sets. While the present construction is also implicitly in terms of cylinder sets (due to Minlos' framework), the reader will notice by comparison that the two probability measures and their support are quite different. As a result the respective kernels take different forms. However, both techniques yield a way to represent the values $h(x)$ for $h$ harmonic and $x \in \verts$ as an integral over ``the boundary''.
 
While Doob's martingale theory works well for harmonic functions in $L^\iy$ or $L^2$, the situation for \HE is different. The primary reason is that \HE is not immediatelly realizable as an $L^2$ space. A considerable advantage of our Gel'fand-Wiener-Ito construction is that \HE is isometrically embedded into $L^2(\Schw',\prob)$ in a particularly nice way: it corresponds to the polynomials of degree 1. See Remark~\ref{rem:Wiener-improves-Minlos}.

Another contrast is that \Lap may, in general, be unbounded in our context. Recall that when studying an operator, an important subtlety is that ``the'' adjoint $\Lap^\ad$ depends on the choice of domain, i.e., the linear subspace $\dom(\Lap) \ci \sH$. %Our domain is markedly different from the usual We consider \Lap as an operator on a rather different Hilbert space, $\ell^2(\verts)$, in \S\ref{sec:L2-theory-of-Lap-and-Trans}.

%There are some sharp contrasts between the case of classical harmonic functions on the disk and  our results for the different classes of harmonic functions on $(\Graph,\cond)$. For instance, in the familiar classical case, the kernel and the geometric boundary (the circle), are identical whether we consider bounded harmonic functions on $D$ or the $L^2$ theory with Hardy space. %Nonetheless, as in the classical theory we will also here obtain limit values for our harmonic functions, and the existence of a limit at boundary points, we show, will exist almost everywhere (a.e.) on $Q = \bd G$.
\end{remark}

Boundary theory is a well-established subject; the deep connections between harmonic analysis, probability, and potential theory have led to several notions of boundary and we will not attempt to give complete references. However, we recommend \cite{Saw97} for introductory material on Martin boundary and \cite{DoSn84, Lyons:ProbOnTrees} for introductory material on resistance networks. Additionally, \cite{TerryLyons, Car73a, Woess00}, and the foundational paper \cite{Nash-Will59} provide more specific background. With regard to infinite graphs and finite-energy functions, see \cite{Soardi94, SoardiWoess91, CaW92, Dod06, PicWoess90, PicWoess88, Wo86, Thomassen90}.
But we ask different questions here, and the operator theory we use is different; it does not easily compare to earlier literature. For some recent related areas, see e.g., \cite{AL08, AAL08, AD06} reproducing kernels, \cite{Arv86} Markov operators, \cite{Cho08} graph analysis, and \cite{AlPa09} operator theory.

There has been a recent interest in analysis and potential theory on infinite-dimensional spaces, and the use of stochastic integration in conjunction with reproducing kernels \cite{HNS09, Xi10, CVTU10,ZXZ09}, and Gel'fand triples \cite{ZX10,BKQ07}. Although our setting here is different, we are able to adapt these tools for the task at hand. This is nontrivial because, in the classical case, there is a natural differentiable structure around, and therefore the choice of a Schwartz space going into a useful Gel'fand triple is often rather conventional.
   But by contrast, we deal with discrete structures, and so we must give up differential operators. Nonetheless, we exhibit Schwartz spaces that yield Gel'fand triples which accomplish what we need.

%%!TEX root = bdG.tex

\section{Basic terms and previous results}
\label{sec:electrical-resistance-networks}

We now proceed to introduce the key notions used throughout this paper: resistance networks, the energy form \energy, the Laplace operator \Lap, and their elementary properties. Our approach is somewhat different from existing studies of networks in the literature, and so we take this opportunity to introduce the tools we will need: an unbounded Laplace operator with dense domain in a Hilbert space, a two-point reproducing kernel for this Hilbert space, the quadratic form associated to the Laplacian, and Gelfand triples. Since these are tools not commonly used in geometric analysis, we include their definitions and some theorems from earlier papers which will be needed later. Additionally, we will use the theorems of Bochner (Theorem~\ref{thm:Bochner's-theorem}), and Minlos (Theorem~\ref{thm:Minlos'-theorem}).

\begin{defn}\label{def:ERN}
  A resistance network is a connected graph $(\Graph,\cond)$, where \Graph is a graph with vertex set \verts, and \cond is the \emph{conductance function} which defines adjacency by $x \nbr y$ iff $c_{xy}>0$, for $x,y \in \verts$. We assume $\cond_{xy} = \cond_{yx} \in [0,\iy)$, and write $\cond(x) := \sum_{y \nbr x} \cond_{xy}$. We require $\cond(x) < \iy$ (note that we allow vertices of infinite degree), but $\cond(x)$ need not be a bounded function on \verts. The notation \cond may be used to indicate the multiplication operator $(\cond v)(x) := \cond(x) v(x)$, 
  \version{}{\marginpar{Is it ok to say ``basis'' here?}}
  i.e., the diagonal matrix with entries $\cond(x)$ with respect to the (vector space) basis $\{\gd_x\}$.
\end{defn}

In this definition, connected means simply that for any $x,y \in \verts$, there is a finite sequence $\{x_i\}_{i=0}^n$ with $x=x_0$, $y=x_n$, and $\cond_{x_{i-1} x_i} > 0$, $i=1,\dots,n$. Conductance is the reciprocal of resistance, so one can think of $(\Graph,\cond)$ as a network of nodes \verts connected by resistors of resistance $\cond_{xy}^{-1}$. We may assume there is at most one edge from $x$ to $y$, as two conductors $\cond^1_{xy}$ and $\cond^2_{xy}$ connected in parallel can be replaced by a single conductor with conductance $\cond_{xy} = \cond^1_{xy} + \cond^2_{xy}$. Also, we assume $\cond_{xx}=0$ so that no vertex has a loop, as electric current will never flow along a conductor connecting a node to itself.\footnote{Nonetheless, self-loops may be useful for technical considerations: one can remove the periodicity of a random walk by allowing self-loops. This can allow one to obtain a ``lazy walk'' which is ergodic, and hence more tractable. See, for example, \cite{LevPerWil08, Lyons:ProbOnTrees}.}

\begin{defn}\label{def:graph-laplacian}
  The \emph{Laplacian} on \Graph is the linear difference operator 
  which acts on a function $v:\verts \to \bR$ by
  \linenopax
  \begin{equation}\label{eqn:def:laplacian}
    (\Lap v)(x) :
    = \sum_{y \nbr x} \cond_{xy}(v(x)-v(y)).
  \end{equation}
  A function $v:\verts \to \bR$ is \emph{harmonic} iff $\Lap v(x)=0$ for each $x \in \verts$. 
\end{defn}

We have adopted the physics convention (so that the spectrum is nonnegative) and thus our Laplacian is the negative of the one commonly found in the PDE literature. The network Laplacian \eqref{eqn:def:laplacian} should not be confused with the stochastically renormalized Laplace operator $\cond^{-1} \Lap$ which appears in the probability literature, or with the spectrally renormalized Laplace operator $\cond^{-1/2} \Lap \cond^{-1/2}$ which appears in the literature on spectral graph theory (e.g., \cite{Chu01}).

\begin{comment}\label{def:Prob-operator}
  The \emph{(probabilistic) transition operator} is defined pointwise for functions on \verts by
  \linenopax
  \begin{align}\label{eqn:def:Prob-trans-oper}%\label{eqn:def:Prob-trans-oper}
    \Prob u(x) = \sum_{y \nbr x} p(x,y) u(y), 
    \q\text{for } 
    p(x,y) = \frac{\cond_{xy}}{\cond(x)},
  \end{align}
  so that $\Lap = \cond(\id - \Prob)$. Note that the harmonic functions are precisely the fixed points of \Prob, and $v=u+k\one$ for $k \in \bC$ implies that $\Prob v = \Prob u + k \one$, so \eqref{eqn:def:Prob-trans-oper} is independent of the representative chosen for $u$.

  The function $p(x,y)$ gives transition probabilities, i.e., the probability that a random walker currently at $x$ will move to $y$ with the next step. Since
  \linenopax
  \begin{align}\label{eqn:def:reversible}
    \cond(x) p(x,y) = \cond(y) p(y,x),
  \end{align}  
  the transition operator \Prob determines a \emph{reversible} Markov process with state space \verts; see \cite{DoSn84,LevPerWil08,Lyons:ProbOnTrees,Peres99}. 
\end{comment}

\begin{defn}\label{def:exhaustion-of-G}
  An \emph{exhaustion} of \Graph is an increasing sequence of finite and connected subgraphs $\{\Graph_k\}_{k=1}^\iy$, so that $\Graph_k \ci \Graph_{k+1}$ and $\Graph = \bigcup \Graph_k$. Since any vertex or edge is eventually contained in some $G_k$, there is no loss of generality in assuming they are contained in $G_1$, for the purposes of a specific computation.
\end{defn}

\begin{defn}\label{def:infinite-vertex-sum}  
  The notation
  \linenopax
  \begin{equation}\label{eqn:def:infinite-sum}
    \sum_{x \in \verts} := \lim_{k \to \iy} \sum_{x \in \Graph_k}
  \end{equation}
  is used whenever the limit is independent of the choice of exhaustion $\{\Graph_k\}$ of \Graph. This is clearly justified, for example, whenever the sum has only finitely many nonzero terms, or is absolutely convergent as in the definition of \energy in Definition~\ref{def:graph-energy}.
\end{defn}

\begin{defn}\label{def:graph-energy}
  The \emph{energy} of functions $u,v:\verts \to \bC$ is given by the (closed, bilinear) Dirichlet form
  \linenopax
  \begin{align}\label{eqn:def:energy-form}
    \energy(u,v)
    := \frac12 \sum_{x \in \verts}  \sum_{y \in \verts} \cond_{xy}(\cj{u}(x)-\cj{u}(y))(v(x)-v(y)),
  \end{align}
  with the energy of $u$ given by $\energy(u) := \energy(u,u)$.
  The \emph{domain of the energy} is
  \linenopax
  \begin{equation}\label{eqn:def:energy-domain}
    \dom \energy = \{u:\verts \to \bC \suth \energy(u)<\iy\}.
  \end{equation}
\end{defn}

Since $\cond_{xy}=\cond_{yx}$ and $\cond_{xy}=0$ for nonadjacent vertices, the initial factor of $\frac12$ in \eqref{eqn:def:energy-form} implies there is exactly one term in the sum for each edge in the network. 

\subsection{The energy space \HE} 
\label{sec:The-energy-space}

 Let \one denote the constant function with value 1 and recall that $\ker \energy = \bC \one$. 

\begin{defn}\label{def:H_energy}\label{def:The-energy-Hilbert-space}
  The energy form \energy is symmetric and positive definite on $\dom \energy$. Then $\dom \energy / \bC \one$ is a vector space with inner product and corresponding norm given by
  \linenopax
  \begin{equation}\label{eqn:energy-inner-product}
    \la u, v \ra_\energy := \energy(u,v)
    \q\text{and}\q
    \|u\|_\energy := \energy(u,u)^{1/2}.
  \end{equation}
  The \emph{energy Hilbert space} \HE is $\dom \energy / \bC \one$ with inner product \eqref{eqn:energy-inner-product}. 
\end{defn}

\begin{defn}\label{def:vx}\label{def:energy-kernel}
  Let $v_x$ be defined to be the unique element of \HE for which
  \linenopax
  \begin{equation}\label{eqn:def:vx}
    \la v_x, u\ra_\energy = u(x)-u(o),
    \qq \text{for every } u \in \HE.
  \end{equation}
  The collection $\{v_x\}_{x \in \verts}$ forms a reproducing kernel for \HE (\cite[Cor.~2.7]{DGG}); we call it the \emph{energy kernel} and \eqref{eqn:def:vx} shows its span is dense in \HE. Note that $v_o$ corresponds to a constant function, since $\la v_o, u\ra_\energy = 0$ for every $u \in \HE$. Therefore, $v_o$ is often ignored or omitted.
\end{defn} 

\begin{defn}\label{def:dipole}
  A \emph{dipole} is any $v \in \HE$ satisfying the pointwise identity $\Lap v = \gd_x - \gd_y$ for some vertices $x,y \in \verts$. One can check that $\Lap v_x = \gd_x - \gd_o$; cf. \cite[Lemma~2.13]{DGG}.
\end{defn}

\begin{defn}\label{def:Fin}
  For $v \in \HE$, one says that $v$ has \emph{finite support} iff there is a finite set $F \ci \verts$ for which $v(x) = k \in \bC$ for all $x \notin F$, i.e., the set of functions of finite support in \HE is 
  \linenopax
  \begin{equation}\label{eqn:span(dx)}
    \spn\{\gd_x\} = \{u \in \dom \energy \suth u(x)=k \text{ for some $k$, for all but finitely many } x \in \verts\},
  \end{equation}
  where $\gd_x$ is the Dirac mass at $x$, i.e., the element of \HE containing the characteristic function of the singleton $\{x\}$. It is immediate from \eqref{eqn:def:energy-form} that $\energy(\gd_x) = \cond(x)$, whence $\gd_x \in \HE$.
  Define \Fin to be the closure of $\spn\{\gd_x\}$ with respect to \energy. 
\end{defn}

\begin{defn}\label{def:Harm}
  The set of harmonic functions of finite energy is denoted
  \linenopax
  \begin{equation}\label{eqn:Harm}
    \Harm := \{v \in \HE \suth \Lap v(x) = 0, \text{ for all } x \in \verts\}.
  \end{equation}
  Note that this is independent of choice of representative for $v$ in virtue of \eqref{eqn:def:laplacian}.
\end{defn}

\begin{lemma}[{\cite[2.11]{DGG}}]
  \label{thm:<delta_x,v>=Lapv(x)}
  For any $x \in \verts$, one has $\la \gd_x, u \ra_\energy = \Lap u(x)$.
\end{lemma}

The following result follows easily from Lemma~\ref{thm:<delta_x,v>=Lapv(x)}; cf.~\cite[Thm.~2.15]{DGG}.

\begin{theorem}[Royden decomposition]
  \label{thm:HE=Fin+Harm}
  $\HE = \Fin \oplus \Harm$.
\end{theorem}

%Theorem~\ref{thm:HE=Fin+Harm} is sometimes called the ``Royden Decomposition'' in honour of the analogous theory established by Royden for Riemann surfaces; see \cite[\S{VI}]{Soardi94}, \cite[\S9.3]{Lyons:ProbOnTrees}. 

\begin{defn}\label{def:monopole}
  A \emph{monopole} at $x \in \verts$ is an element $w_x \in \HE$ which satisfies
  $\Lap w_x(y) = \gd_{xy}$, where $\gd_{xy}$ is Kronecker's delta. 
  %By \cite[Lemma~3.2]{DGG}, this is equivalent to $v = \LapV^\ad u$ in \HE, that is,
  %\linenopax
  %\begin{equation}\label{eqn:def:monopole}
  %  \la w_x, \Lap u \ra_\energy = \la \gd_x, u \ra_\energy,
  %  \qq \text{ for all } u \in \dom \LapV.
  %\end{equation}
  %When nonempty, the space of monopoles at the origin is closed and convex, so \energy attains a unique minimum here. 
  In case the network supports monopoles, let $w_o$ always denote the unique energy-minimizing monopole at the origin. %, and denote $\monov := v_x + w_o$ and $\monof := f_x + w_o$, where $f_x = \Pfin v_x$. %By \cite[Cor.~3.7]{DGG}, one has $\monov = \monof$ for all $x$ iff if $\Harm=0$. 
 %Note that $w_o \in \Fin$.
  
  When \HE contains monopoles, let $\MP_x$ denote the vector space spanned by the monopoles at $x$. This implies that $\MP_x$ may contain harmonic functions; see \cite[Lemma~4.1]{DGG}. With $v_x$ and $f_x = \Pfin v_x$ as above, we indicate the distinguished monopoles 
  \linenopax
  \begin{align}\label{eqn:def:monov-and-monof}
    \monov := v_x + w_o
    \q\text{and}\q
    \monof := f_x + w_o.
  \end{align}
  % (Corollary~\ref{thm:Harm-nonzero-iff-multiple-monopoles} below confirms that $\monov = \monof$ for all $x$ iff if $\Harm=0$.) 
\end{defn}

\begin{remark}\label{rem:w_O-in-Fin}
  Note that $w_o \in \Fin$, whenever it is present in \HE, and similarly that \monof is the energy-minimizing element of $\MP_x$. To see this, suppose $w_x$ is any monopole at $x$. Since $w_x \in \HE$, write $w_x = f+h$ by Theorem~\ref{thm:HE=Fin+Harm}, and get $\energy(w_x) = \energy(f) + \energy(h)$. Projecting away the harmonic component will not affect the monopole property, so $\monof = \Pfin w_x$ is the unique monopole of minimal energy.  The Green function is $g(x,y) = w_y^o(x)$, where $w_y^o$ is the representative of \monof[y] satisfying $w_y^o(o)=0$. %Also, $w_o$ corresponds to the projection of \one to \Gddo; see \S\ref{sec:grounded-energy-space}.
\end{remark}

\begin{defn}\label{def:LapM}
  The dense subspace of \HE spanned by monopoles (or dipoles) is
  \linenopax
  \begin{equation}\label{eqn:def:MP}
    \MP := \spn\{v_x\}_{x \in \verts} + \spn\{\monov, \monof\}_{x \in \verts}.
  \end{equation}
 % When $\MP_x \neq \es$, and we say that $u \in \MP$ if and only if $u \in \dom \LapM^\ad$ and $u$ is a (finite) linear combination of functions $w_{x_1}, \dots, w_{x_n}$, each satisfying $\Lap w_{x_i} = \gd_{x_i}$. 
  %Note that for an element $w = \sum_{i=1}^n a_i w_{x_i} \in \MP$, $w_{x_i}$ need not have finite energy. As a consequence, \MP always contains the energy kernel $\{v_x\}$, even when there are no monopoles in \HE. In fact, it is easy to see that if there are no monopoles in \HE (i.e., if the network is recurrent), then $\MP = \spn \{v_x\}$. 
   
 %Let $\sV := \spn\{v_x\}_{x \in \verts}$ denote the vector space of \emph{finite} linear combinations of dipoles from the energy kernel, and 
  Let \LapM be the closure of the Laplacian when taken to have the dense domain \MP. Since \Lap agrees with \LapM pointwise, we may suppress reference to the domain for ease of notation.
\end{defn}

\begin{lemma}[{\cite[Lemma~3.5]{DGG}}]
  \label{thm:LapM-is-semibounded}
  \LapM is Hermitian with $\la u, \LapM u\ra_\energy \geq 0$ for all $u \in \MP$.
\end{lemma}

\begin{lemma}[{\cite[Lemma~3.6]{DGG}}]
  \label{thm:MP-contains-spans}
  When the network is transient, \MP contains the spaces $\spn\{v_x\}, \spn\{f_x\}$, and $\spn\{h_x\}$, where $f_x = \Pfin v_x$ and $h_x = \Phar v_x$. When the network is not transient, $\MP = \spn\{v_x\}= \spn\{f_x\}$.
\end{lemma}

\begin{remark}[Monopoles and transience]
  \label{rem:transient-iff-monopoles}
  The presence of monopoles in \HE is equivalent to the transience of the simple random walk on the network with transition probabilities $p(x,y) = \cond_{xy}/\cond(x)$: note that if $w_x$ is a monopole, then the current induced by $w_x$ is a unit flow to infinity with finite energy. It was proved in \cite{TerryLyons} that the network is transient if and only if there exists a unit current flow to infinity; see also \cite[Thm.~2.10]{Lyons:ProbOnTrees}. %It is also clear that the existence of a monopole at one vertex is equivalent to the existence of a monopole at every vertex: consider $v_x + w_o$. The corresponding statement about transience is well-known.
\end{remark}

\subsection{The discrete Gauss-Green identity}
\label{sec:The-discrete-Gauss-Green-identity}

The space \MP is introduced as a dense domain for \Lap and as the scope of validity for the discrete Gauss-Green identity of Theorem~\ref{thm:Discrete-Gauss-Green-identity}. 

\begin{defn}\label{def:subgraph-boundary}
  If $H$ is a subgraph of $G$, then the boundary of $H$ is
  \linenopax
  \begin{equation}\label{eqn:subgraph-boundary}
    \bd H := \{x \in H \suth \exists y \in H^\complm, y \nbr x\}.
  \end{equation}
  \glossary{name={$\bd H$},description={boundary of a subgraph},sort=b,format=textbf}
  The \emph{interior} of a subgraph $H$ consists of the vertices in $H$ whose neighbours also lie in $H$:
  \linenopax
  \begin{equation}\label{eqn:interior}
    \inn H := \{x \in H \suth y \nbr x \implies y \in H\} = H \less \bd H.
  \end{equation}
  \glossary{name={$\inn H$},description={interior of a subgraph},sort=i,format=textbf}
  For vertices in the boundary of a subgraph, the \emph{normal derivative} of $v$ is
  \linenopax
  \begin{equation}\label{eqn:sum-of-normal-derivs}
    \dn v(x) := \sum_{y \in H} \cond_{xy} (v(x) - v(y)),
    \qq \text{for } x \in \bd H.
  \end{equation}
  Thus, the normal derivative of $v$ is computed like $\Lap v(x)$, except that the sum extends only over the neighbours of $x$ which lie in $H$.
\end{defn}

Definition~\ref{def:subgraph-boundary} will be used primarily for subgraphs that form an exhaustion of \Graph, in the sense of Definition~\ref{def:exhaustion-of-G}: an increasing sequence of finite and connected subgraphs $\{\Graph_k\}$, so that $\Graph_k \ci \Graph_{k+1}$ and $\Graph = \bigcup \Graph_k$. Also, recall that $\sum_{\bd \Graph} := \lim_{k \to \iy} \sum_{\bd \Graph_k}$ from Definition~\ref{def:boundary-sum}.

\begin{defn}\label{def:boundary-sum}
  A \emph{boundary sum} is computed in terms of an exhaustion $\{G_k\}$ by
  \linenopax
  \begin{equation}\label{eqn:boundary-sum}
    \sum_{\bd \Graph} := \lim_{k \to \iy} \sum_{\bd \Graph_k},
  \end{equation}
  whenever the limit is independent of the choice of exhaustion, as in Definition~\ref{def:infinite-vertex-sum}.
\end{defn}

On a finite network, all harmonic functions of finite energy are constant, so that $\HE = \Fin$ by Theorem~\ref{thm:HE=Fin+Harm}, and one has $\energy(u,v) = \sum_{x \in \verts} u(x) \Lap v(x)$, for all $u,v \in \HE$. In fact, this remains true for recurrent infinite networks, as shown in \cite[Thm.~4.4]{DGG}; see also \cite{Kayano88}. However, the possibilities are much richer on an infinite network, as evinced by the following theorem.

\begin{theorem}[Discrete Gauss-Green identity]
  \label{thm:Discrete-Gauss-Green-identity}
  If $u \in \HE$ and $v \in \MP$, then
  \linenopax
  \begin{equation}\label{eqn:E(u,v)=<u_0,Lapv>+sum(normals)}
    \la u, v \ra_\energy
    = \sum_{\verts} \cj{u} \Lap v
      + \sum_{\bd \Graph} \cj{u} \dn v.
  \end{equation}
\end{theorem}

The discrete Gauss-Green formula \eqref{eqn:E(u,v)=<u_0,Lapv>+sum(normals)} is the main result of \cite{DGG}; that paper contains several consequences of the formula, especially as pertains to transience.

%from PALLE:
%Here we are motivated by data located on infinite graphs and its analysis, while in physics the goal is to represent systems with an infinite number of particles, or quantum fields; second quantization. The two approaches have common aims: the use of probability measures P in solving equations; but the probability measures as it turns out must be carried by certain dual spaces, larger than the initial Hilbert space. The technical terms are nuclear spaces $S$ embedded in Hilbert space: Now $S$ will have its own topology, but it will be a dense linear subspace in $H$ (the energy space \HE here!). As a result the dual $SÕ$ of $S$ will contain the Hilbert space $H$. The three spaces together are called Gelfand-triples. In the better known instances of this, $S$ will be a Schwartz space of test functions, $H$ an $L^2$ space, and $SÕ$ will be a space of tempered distributions. As it turns out our boundaries $\bd G$ may be realized in $SÕ$ for such a Gelfand construction.

\subsection{Gel'fand triples and duality}
\label{sec:Gel'fand-triples-and-duality}

One would like to obtain a (probability) measure space to serve as the boundary of \Graph. It is shown in \cite{Gro67,Gro70,Minlos63} that no Hilbert space of functions \sH is sufficient to support a Gaussian measure \prob (i.e., it is not possible to have $0<\prob(\sH)<\iy$ for a \gs-finite measure). However, it \emph{is} possible to construct a \emph{Gel'fand triple} (also called a \emph{rigged Hilbert space}): a dense subspace $S$ of \sH with
\linenopax
\begin{equation}\label{eqn:Gel'fand-triple-intro}
  S \ci \sH \ci S',
\end{equation}
where $S$ is dense in \sH and $S'$ is the dual of $S$. Additionally, $S$ and $S'$ must also satisfy some technical conditions: $S$ is a 
%Fr\'{e}chet space in its own right but realized as 
dense subspace of \sH with respect to the Hilbert norm, but also comes equipped with a strictly finer ``test function'' topology. When $S$ is a Fr\'echet space equipped with a countable system of seminorms (stronger than the norm on \sH), then the inclusion map of $S$ into \sH is continuous; in fact, it is possible to chose the seminorms in such a way that one gets a nuclear embedding (details below). 
%It is assumed that the inclusion mapping of $S$ into \sH is continuous in the respective topologies. 
Therefore, when the dual $S'$ is taken with respect to this finer (Fr\'echet) topology, one obtains a strict containment $\sH \subsetneq S'$. It turns out that $S'$ is large enough to support a nice probability measure, even though \sH is not.

It was Gel'fand's idea to formalize this construction abstractly using a system of nuclearity axioms \cite{GMS58, Minlos58, Minlos59}. Our presentation here is adapted from quantum mechanics and the goal is to realize $\bd \Graph$ as a subset of $S'$. We will give a ``test function topology'' as a Fr\'{e}chet topology defined via a specific sequence of seminorms, using either an onb for \HE (in \S\ref{sec:Gel'fand-triple-for-HE}) or the domain of $\Lap^p$ (in \S\ref{sec:unbounded-case}).

\begin{remark}[Tempered distributions and the Laplacian]
  \label{rem:tempered-distributions}
  There is a concrete situation when the Gel'fand triple construction is especially natural: $\sH = L^2(\bR,dx)$ and $S$ is the \emph{Schwartz space} of functions of rapid decay. That is, each $f \in S$ is $C^\iy$ smooth functions which decays (along with all its derivatives) faster than any polynomial. In this case, $SÕ$ is the space of \emph{tempered distributions} and the seminorms defining the Fr\'{e}chet topology on $S$ are 
\linenopax
\begin{align*}
  p_m(f) := \sup \{|x^k f^{(n)}(x)| \suth x \in \bR, 0 \leq k,n \leq m\},
  \qq m=0,1,2,\dots,
\end{align*}
where $f^{(n)}$ is the \nth derivative of $f$. Then $S'$ is the dual of $S$ with respect to this Fr\'{e}chet topology. One can equivalently express $S$ as
\linenopax
\begin{align}\label{eqn:Schwartz-space-as-powers-of-Hamiltonian}
  S := \{f \in L^2(\bR) \suth (\tilde{P}^2 + \tilde{Q}^2)^n f \in L^2(\bR), \forall n\},
\end{align}
where $\tilde{P}:f(x) \mapsto \frac1\ii\frac{d}{dx}$ and $\tilde{Q}:f(x) \mapsto x f(x)$ are Heisenberg's operators. The operator $\tilde{P}^2 + \tilde{Q}^2$ is often called the quantum mechanical Hamiltonian, but some others (e.g., Hida, Gross) would call it a Laplacian, and this perspective tightens the analogy with the present study. In this sense, \eqref{eqn:Schwartz-space-as-powers-of-Hamiltonian} could be rewritten $S := \dom \Lap^\iy$; compare to \eqref{eqn:Schwartz}.
We show that a general network $(\Graph,\cond)$ always has a harmonic oscillator; in fact, we discuss an operator \Onb in Definition~\ref{def:Onb} which is unitarily equivalent to $\tilde{P}^2 + \tilde{Q}^2$ and hence has the same spectrum. %) with key properties in common with the better known variant from Heisenberg's quantum mechanics.
\end{remark}

  The duality between $S$ and $S'$ allows for the extension of the inner product on \sH to a pairing of $S$ and $S'$:
\linenopax
\begin{align*}
  \la \cdot,\cdot\ra_\sH:\sH \times \sH \to \bC
  \qq\text{to}\qq
  \la \cdot,\cdot\ra_\sH^\sim: S \times S' \to \bR.
\end{align*}
In other words, one obtains a Fourier-type duality restricted to $S$. Moreover, it is possible to construct a Gel'fand triple in such a way that $\prob(S')=1$ for a Gaussian probability measure \prob. When applied to $\sH=\HE$, the construction yields three main outcomes:
\begin{enumerate}
  \item The next best thing to a Fourier transform for an arbitrary graph.
  \item A concrete representation of \HE as an $L^2$ measure space $\HE \cong L^2(S',\prob)$.
  \item A boundary integral representation for harmonic functions of finite energy.
\end{enumerate}
%\version{}{This offers the exciting possibility for Fourier analysis of nonabelian groups whose Cayley graphs conform to our definition of an \ERN. The authors are currently plumbing the viability of this approach in another paper.}

%The reader may find much useful background information for this material in \cite{Arv76a,Arv76b}. 
As a prelude, we begin with Bochner's Theorem, which characterizes the Fourier transform of a positive finite Borel measure on the real line. The reader may find \cite{ReedSimonII} helpful for further information.
\begin{theorem}[Bochner]
  \label{thm:Bochner's-theorem}
  Let $G$ be a locally compact abelian group. Then there is a bijective correspondence $\sF:\sM(G) \to \sP\sD(\hat G)$, where $\sM(G)$ is the collection of measures on $G$, and $\sP\sD(\hat G)$ is the set of positive definite functions on the dual group of $G$. Moreover, this bijection is given by the Fourier transform
  \linenopax
  \begin{equation}\label{eqn:Bochner's-theorem}
    \sF:\gn \mapsto \gf_\gn
    \qq\text{by}\qq
    \gf_\gn(\gx) = \int_G e^{\ii \la \gx,x\ra} \,d\gn(x).
  \end{equation}
\end{theorem}
%In our applications to the \ERN $(\bZd,\one)$ in \S\ref{sec:lattice-networks}, the underlying group structure allows us to apply the above version of Bochner's theorem. Specifically, in the context of group duality, Bochner's theorem characterizes the Fourier transform of a positive finite Borel measures; cf.~\cite{ReedSimonII,Ber96}.

For our representation of the energy Hilbert space \HE in the case of general resistance network, we will need Minlos' generalization of Bochner's theorem from \cite{Minlos63, Sch73}. This important result states that a cylindrical measure on the dual of a nuclear space is a Radon measure iff its Fourier transform is continuous. In this context, however, the notion of Fourier transform is infinite-dimensional, and is dealt with by the introduction of Gel'fand triples \cite{Lee96}.

\begin{theorem}[Minlos]
  \label{thm:Minlos'-theorem}
  Given a Gel'fand triple $S \ci \sH \ci S'$, Bochner's Theorem may be extended to yield a bijective correspondence between the positive definite functions on $S$ and the Radon probability measures on $S'$. Moreover, in a specific case, this correspondence is uniquely determined by the identity
  \linenopax
  \begin{align}\label{eqn:Minlos-identity}
    \int_{S'} e^{\ii \la u, \gx\ra_{\tilde \sH}} \,d\prob(\gx) 
    = e^{-\frac12 \la u,u\ra_\sH},
  \end{align}
  where $\la \cdot , \cdot \ra_\sH$ is the original inner product on \sH and $\la \cdot , \cdot \ra_{\tilde \sH}$ is its extension to the pairing on $S \times S'$. 
\end{theorem}

Formula \eqref{eqn:Minlos-identity} may be interpreted as defining the Fourier transform of \prob; the function on the right-hand side is positive definite and plays a special role in stochastic integration, and its use in quantization. 
%%!TEX root = bdG.tex

\section{Gel'fand triples for \HE}
\label{sec:Gel'fand-triple-for-HE}

In this section, we describe two methods for construction a Gel'fand triple for \HE. The first method is applicable to all networks, but relies on the choice of some enumeration of the vertices of \Graph, and the Gram-Schmidt algorithm for producing an onb. However, we will see that the Gram-Schmidt algorithm yields a much more explicit formula than usual, in the present context. The second method is applicable only when the Laplacian is unbounded. However, in this case the construction does not require any enumeration (or onb) and may provide for more feasible computations.

\begin{remark}\label{rem:Schwartz-space-is-real-valued}
  Note that \Schw and $\Schw'$ consist of \bR-valued functions in this section. This technical detail is important because we do not expect the integral $\int_{S'} e^{\ii \la u, \cdot \ra_{\tilde \Wiener}} \,d\prob$ from \eqref{eqn:Minlos-identity} to converge unless it is certain that $\la u, \cdot \ra$ is \bR-valued. %This is the reason for the last conclusion of Lemma~\ref{thm:energy-extends}.
  After the Wiener embedding is carried out in Theorem~\ref{thm:HE-isom-to-L2(S',P)}, all results can be complexified.
\end{remark}

\subsection{Gel'fand triples via Gram-Schmidt}
\label{sec:Gel'fand-triples-via-Gram-Schmidt}

In this section, we describe a Gel'fand triple for \HE where the class of test functions \Schw is described in terms of the decay properties of a certain orthonormal basis (onb) for \HE. We will see in Remark~\ref{rem:onb-gives-Gaussian-iid} that this onb corresponds to a system of i.i.d. random variables (which are, in fact, Gaussian with mean $0$ and variance $1$).

The onb $\{\onb_n\}_{n \in \bN}$ comes by applying the Gram-Schmidt process to the reproducing kernel $\{v_{x_n}\}_{n \in \bN}$, where we have fixed some enumeration $\{x_n\}_{n \in \bN}$ of the vertices $\Graph\less\{o\}$. That is, we put $x_0=o$ and henceforth exclude $x_0$ from the discussion, as it will not be relevant.
Given $\{\onb_1,\dots,\onb_{n-1}\}$, one obtains $\onb_n$ via
\linenopax
\begin{align}\label{eqn:Gram-Schmidt-matrix}
  \left[\begin{array}{c}
    \onb_1 \\ \onb_2 \\ \onb_3 \\ \vdots \\ \onb_n
  \end{array}\right]
  =
  \left[\begin{array}{cccccc}
    \|v_{x_1}\|_\energy^{-1} & 0 & 0 & 0 & \dots & 0 \\ 
    M_{2,1} & M_{2,2} & 0 & 0 & \dots & 0 \\ 
    M_{3,1} & M_{3,2} & M_{3,3} & 0 & \dots & 0 \\ 
    \vdots \\ 
    M_{n,1} & M_{n,2} & M_{n,3} & \dots & \dots & M_{n,n}
  \end{array}\right]
  \left[\begin{array}{c}
    v_{x_1} \\ v_{x_2} \\ v_{x_3} \\ \vdots \\ v_{x_n}
  \end{array}\right].
\end{align}
Consequently, for each $N \in \bN$, the triangular nature of $M$ gives
\linenopax
\begin{align}\label{eqn:dipole-span=basis-span}
  \spn\{v_{x_1},\dots,v_{x_N}\} = \spn\{\onb_1,\dots,\onb_N\}.
\end{align}

\begin{remark}\label{rem:Gram-inverse-formula}
  Note that the reproducing kernel gives one an explicit formula for the entries of the inverse of this particular Gram-Schmidt matrix: 
  \linenopax
  \begin{align}\label{eqn:M-inv-entries}
    (M^{-1})_{i,j} = \la v_{x_i},\onb_j\ra_\energy = \onb_j(x_i) - \onb_j(o). 
  \end{align}
  This is certainly in distinct contrast with the general case, and allows us to find a formula for the entries of $M$ itself in Lemma~\ref{thm:GS-entries}.
\end{remark}

%\begin{defn}\label{def:x<=y}
%  We use the notation $x \leq y$ to indicate that $x$ occurs not after $y$ in the (fixed) enumeration of $\Graph \less \{o\}$. More precisely, one has $x=x_n$ and $y=y_m$ with $n \leq m$.  
%\end{defn}

\begin{lemma}\label{thm:GS-entries}
  The entries of the Gram-Schmidt matrix $M$ are given by
  \linenopax
  \begin{align}\label{eqn:GS-entries}
    M_{i,j} = 
    \begin{cases}
      (\Lap \onb_i)(x_j), & j \leq i \\
      0 & \text{else},
    \end{cases}
    \qq\text{for } i,j=1,2,\dots. 
  \end{align}
  %where $y \leq x$ is as in Definition~\ref{def:x<=y}.
  \begin{proof}
    For $j \leq i$, an application of \eqref{thm:<delta_x,v>=Lapv(x)} gives
    \linenopax
    \begin{align*}%\label{eqn:GS-entries-computation}
      \Lap \onb_i(x_j)
      &= \la \gd_{x_j},\onb_i\ra_\energy 
      = \left\la \gd_{x_j}, \sum_{k \leq i} M_{i,k} v_{x_k} \right \ra_\energy 
      = \sum_{k \leq i} M_{i,k} \la \gd_{x_j}, v_{x_k}\ra_\energy 
      = \sum_{k \leq i} M_{i,k} (\gd_{x_j}(x_k)-\gd_{x_j}(o)), 
    \end{align*}
    where the last equality comes by \eqref{eqn:def:vx}.
    Note that $\gd_x(y)=\gd_x(o)$ for every $y$ except $y=x$, 
    %$M$ does not contain a row (or column) corresponding to $o$, as noted just before \eqref{eqn:Gram-Schmidt-matrix}, 
    so the last sum above has a nonzero term only for $k=j$, and the result follows.
  \end{proof}
\end{lemma}

\begin{comment}\label{rem:M-inv-product}
  The reader can verify that $(MM^{-1})_{x,y} = (M^{-1}M)_{x,y} = \gd_{x,y}$ (Kronecker \gd) as follows:
  \linenopax
  \begin{align}\label{eqn:M-inv-product}
    \sum_{z \leq x \wedge y} \Lap \onb_z(x) \onb_z(y) 
    = \sum_{z \leq x \wedge y} \Lap (\onb_z(x) -\onb_z(o))\onb_z(y)
    = \Lap\left(\sum_{z \leq x \wedge y} \onb_z(y) (\onb_z -\onb_z(o))\right)(x)
  \end{align}
  where the last is Kronecker's \gd.
\end{comment}

From \eqref{eqn:M-inv-entries} and Lemma~\ref{thm:GS-entries}, we have the handy conversion formulas: 
\linenopax
\begin{align}\label{eqn:conversion-formulas}
  \onb_i = \sum_{j \leq i} \Lap \onb_i(x_j) v_{x_j}
  \qq\text{and}\qq
  v_{x_i} = \sum_{j \leq i}  \left(\onb_k(x_i)\vstr[1.7] - \onb_k(o)\right) \onb_j.
\end{align}

\begin{lemma}\label{thm:Kronecker-sum}
  We have the identity
  \linenopax
  \begin{align}\label{eqn:Kronecker-sum}
    \sum_{j \leq k \leq i} \left(\onb_k(x_i)\vstr[1.7] - \onb_k(o)\right) \Lap \onb_k (x_j)
    = \gd_{i,j},
    \qq \text{for }i,j=1,2,\dots.
  \end{align}
  \begin{proof}
    By formula \eqref{eqn:M-inv-entries}, the left side of \eqref{eqn:Kronecker-sum} is equal to
    \linenopax
    \begin{align*}%\label{eqn:}
      \sum_{j \leq k \leq i} \left(\onb_k(x_i)\vstr[1.7] - \onb_k(o)\right) \Lap \onb_k (x_j)
      &= \Lap \left( \sum_{k \leq i} \la v_{x_i},\onb_j\ra_\energy \onb_k (x_j)\right) 
      = \Lap v_{x_i}(x_j) 
      = \gd_{x_i}(x_j) - \gd_{o}(x_j).
    \end{align*}
    Note that $\Lap \onb_k(x_j)=0$ for $j>k$, so the second sum runs over all $k \leq i$. Also, note that $\gd_{x_i}(x_j) - \gd_{o}(x_j) = \gd_{i,j}$ for $i,j>0$ (and the indexing of $M$ begins at $1$, not $0$).
  \end{proof}
\end{lemma}
Lemma~\ref{thm:Kronecker-sum} can also be proven by combining the identities in \eqref{eqn:conversion-formulas}.

\begin{lemma}\label{thm:V=EE*}
  Let $V_{x,y} := \la v_x, v_y\ra_\energy$, and let $E = M^{-1}$ be defined as in \eqref{eqn:M-inv-entries}. Then $EE^\ad = V$.
  \begin{proof}
    Computing entrywise,
    \begin{align*}%\label{eqn:}
      (EE^\ad)_{i,j}
      = \sum_{k} E_{x_i,x_k} E_{x_j,x_k} 
      = \sum_{k} \left(\onb_k(x_i)\vstr[1.8] - \onb_k(o)\right)\left(\onb_k(x_j)\vstr[1.7] - \onb_k(o)\right) 
      &= \sum_{k} \la v_{x_i},\onb_k\ra_\energy \la v_{x_j},\onb_k\ra_\energy,
      %&= \la v_{x_i}, v_{_j}\ra_\energy
      %= V_{x_i,x_j}.
    \end{align*}
    which is equal to $\la v_{x_i}, v_{x_j}\ra_\energy$ by Parseval's identity.
  \end{proof}
\end{lemma}

\begin{defn}\label{def:S,S'}
  The space of test functions and the space of distributions corresponding to the onb $\{\onb_n\}_{n\in \bN}$ are defined by
  \linenopax
  \begin{align}
    \Schw &= \{s = \sum_{n \in \bN} s_n \onb_n \suth \forall p \in \bN, \exists C > 0 \text{ such that } |s_n| \leq C/n^p\}, \text{ and} \label{eqn:def:S} \\
    \Schw' &= \{\gx = \sum_{n \in \bN} \gx_n \onb_n \suth \exists p \in \bN, \exists C > 0 \text{ such that } |\gx_n| \leq Cn^p\} \label{eqn:def:S'}.
  \end{align} 
  Thus, $\Schw = \bigcap_{p \in \bN} \{s \suth \|s\|_p < \iy\}$ where the Fr\'{e}chet \emph{$p$-seminorm} of $s = \sum_{n \in \bN} s_n \onb_n$ is 
  \linenopax
  \begin{align}\label{eqn:p-seminorm}
     \|s\|_p 
     := \left(\sum_{n \in \bN} n^p |s_n|^2\right)^{1/2},
      \qq s \in \Schw, p \in \bN.
  \end{align}
  Note that the system of seminorms \eqref{eqn:p-seminorm} is equivalent to system of seminorms defined by 
  \linenopax
  \begin{align}\label{eqn:p-seminorm2}
     \|s\|_{p} := \sup_{n \in \bN} n^p |s_n|,
      \qq s \in \Schw, p \in \bN,
  \end{align}
  in the sense that both define the same Fr\'echet topology on \Schw. (Each seminorm in one system is dominated by one from the other, but with a different $p$.) We occasionally find it more convenient to calculate with \eqref{eqn:p-seminorm2} instead of \eqref{eqn:p-seminorm}.
\end{defn}

\begin{defn}\label{def:Onb}
  Let $\sV := \spn\{v_x\}_{x \in \Graph}$ and define a mapping $\Onb:\sV \to \HE$ by 
  \linenopax
  \begin{align}\label{eqn:Onb}
    \Onb v_{x_n} = \sum_{k=1}^n k \onb_k(x_n) \onb_k.
  \end{align}
\end{defn}

\begin{remark}\label{rem:Onb-formulas}
  From \eqref{eqn:Onb}, one has 
  \linenopax
  \begin{align}\label{eqn:|Onbv|}
    \|\Onb v_{x_n}\|_\energy^2 = \sum_{k=1}^n k^2 |\onb_k(x_n)|^2
    \qq\text{and}\qq
    \la v_{x_n}, \Onb v_{x_m}\ra_\energy 
      = \sum_{k=1}^{n \wedge m} k \onb_k(x_n) \onb_k(x_m).
  \end{align}
  Note that $\onb_k \in \sV$ by \eqref{eqn:dipole-span=basis-span}, and that $\Onb \onb_k = k \onb_k$ for each $k \in \bN$. We use the symbol \Onb for the operator discussed in this section by way of analogy with the number operator N from quantum mechanics. Indeed, \Onb can also be defined as $a^\ast a$ for a suitable operator $a$ and its adjoint.
\end{remark}

In the following lemma, we use the symbol $\bar{\Onb}$ to denote the closure of the operator \Onb (i.e., the domain is the closure of $\spn\{\onb_n\}$ with respect to the graph norm). 

\begin{lemma}\label{thm:properties-of-Onb}
  The mapping \Onb is essentially self-adjoint, and is unbounded if and only if \Graph is infinite. Moreover, if we define the seminorms $\gr_n(u) := \|(\bar{\Onb})^nu\|_\energy$, then $\{\gr_n\}$ and $\{\|\cdot\|_p\}$ induce equivalent topologies on \Schw, so that
  \linenopax
  \begin{align}\label{eqn:Schw-from-Onb}
    \Schw = \bigcap_{n \in \bN} \dom (\bar{\Onb})^n 
  \end{align}
  and $u \in \Schw$ if and only if $\gr_n(u) < \iy$ for each $n \in \bN$.
  \begin{proof}
    Unitary equivalence of \HE with $\ell^2(\bZ_+)$ is given by $U:\onb_n \mapsto \gd_n$, where $\gd_n(m) := \gd_{n,m}$ (Kronecker \gd) for $n,m \in \bZ_+$. Define $\bL:\ell^2(\bZ_+) \to \ell^2(\bZ_+)$ by $\bL \gd_n = n \gd_n$ so that $\bL U = U\Onb$ holds on the dense subspace $\spn\{\onb_n\}$. The rest follows by \cite[\S2]{Simon}.
  \end{proof}
\end{lemma}

\begin{cor}
  \label{thm:sandwich=triple}
  The inclusion mapping $\Schw \hookrightarrow \HE$ is nuclear, and so $\Schw \ci \HE \ci \Schw'$ is a Gel'fand triple.
  \begin{proof}
    When the space of test functions is defined as $\dom T^\iy$ for some operator $T$ with pure point spectrum (as in \eqref{eqn:Schw-from-Onb}), then nuclearity follows if there is a $p \in \bZ_+$ such that the reciprocal eigenvalues of $T$ are $p$-summable; see \cite{Simon}. Since \Onb has spectrum $\bZ_+$ and $\sum_{n=1}^\iy n^{-p} < \iy$ for $p \geq 2$, the conclusion follows.
  \end{proof}
\end{cor}

\begin{lemma}\label{thm:v_x-in-S}
  The energy kernel $\{v_x\}_{x \in \verts}$ is a Fr\'{e}chet-dense subset of \Schw. 
  \begin{proof}
    In the expansion with respect to the onb as in \eqref{eqn:def:S}, the basis element $\onb_k$ has coefficients $s_n = \gd_{nk}$. Since this sequence $\{s_n\}$ vanishes after $n=k$, it clearly satisfies the required decay condition $|s_n| \leq Cn^{-p}$. From \eqref{eqn:dipole-span=basis-span}, the same clearly holds for $v_{x_k}$. This shows that the kernel is contained in \Schw.
    
    To see that $\{v_x\}$ is dense in \Schw, it suffices by \eqref{eqn:dipole-span=basis-span} to show that the onb $\{\onb_n\}_{n \in \bN}$ is dense. 
    %Fix any $p \in \bN$ and $\ge > 0$. 
    Given any $u = \sum u_k \ge_k \in \Schw$ and any $p \in \bN$,
    %we must find an $f$ for which $\|u-f\|_p < \ge$. Since $u \in \Schw$, 
    there is a $C = C_p$ such that $|u_k| \leq C/k^{p+1}$. Now if $u_N = \sum_{k=1}^N u_k \onb_k$ is the \nth[N] truncation of $u$, then
    \linenopax
    \begin{align*}%\label{eqn:}
      \|u-f\|_p
      = \sup_{k} k^p |u_k - f_k|
      = \sup_{k > N} k^p |u_k - f_k|
      \leq k^p \frac{C}{k^{p+1}} \limas{k} 0.
    \end{align*}
    Thus, one can always approximate $u \in \Schw$ by $u_N \in \spn\{\onb_1,\dots,\onb_N\}$. % with respect to the Fr\"{e}chet topology.
  \end{proof}
\end{lemma}

\begin{comment}\label{thm:T-injective}
  Let $T:\prod_{\Schw} \bR \to \prod_{\verts} \bR$ be the transformation defined by $T(\gx)(x) = \la v_x, \gx \ra$ for $\gx \in \Schw'$. Then $T$ is injective.
  \begin{proof}
    Let $\gx, \gz \in \Schw'$ and suppose $\la v_x, \gx \ra_\energy = \la v_x, \gz \ra_\energy$ for all $x$. Then $\gx(v_x) = \gz(v_x)$ implies $\gx=\gz$ on a Fr\'echet-dense subset of \Schw, and hence on all of \Schw by continuity.
  \end{proof}
\end{comment}

\subsection{Gel'fand triples in the case when \Lap is unbounded}
\label{sec:unbounded-case}

In the case when $\Lap:\HE \to \HE$ is unbounded, there is an alternative construction of $\Schw$ and $\Schw'$, which begins by identifying a certain subspace of $\MP = \dom \LapM$ (as given in Definition~\ref{def:LapM}) to act as the space of test functions.

\begin{defn}\label{def:extn-of-Lap}
  Let \LapS be a self-adjoint extension of \LapM; since \LapM is Hermitian and commutes with conjugation (since \cond is \bR-valued), a theorem of von Neumann's states that such an extension exists.
\glossary{name={\LapS},description={a self-adjoint extension of \LapM},sort=L,format=textbf}
  
  Let $\LapS^p u := (\LapS\LapS\dots\LapS)u$ be the $p$-fold product of \LapS applied to $u \in \HE$. Define $\dom(\LapS^p)$ inductively by
  \linenopax
  \begin{equation}\label{eqn:domE(LapEp)}
    \dom(\LapS^p) := \{u \suth \LapS^{p-1} u \in \dom(\LapS)\}.
  \end{equation}
\end{defn}

\begin{defn}\label{def:Schw-on-G}
  The \emph{(Schwartz) space of potentials of rapid decay} is 
  \linenopax
  \begin{equation}\label{eqn:Schwartz}
    \Schw := \dom(\LapS^\iy), 
  \end{equation}
\glossary{name={$\Schw$},description={``Schwartz space'' of test functions (of rapid decay)},sort=S,format=textbf}
  where $\dom(\LapS^\iy) := \bigcap_{p=1}^\iy \dom(\LapS^p)$ consists of all \bR-valued functions $u \in \HE$ for which $\LapS^p u \in \HE$ for any $p$. The space of \emph{Schwartz distributions} or \emph{tempered distributions} is the dual space $\Schw'$ of \bR-valued continuous linear functionals on \Schw.
\end{defn}

\begin{remark}\label{rem:Schw-is-Frechet}
  Note that \Schw is dense in $\dom(\LapS)$ with respect to the graph norm, by standard spectral theory. For each $p \in \bN$, there is a seminorm on \Schw defined by
  \linenopax
  \begin{equation}\label{eqn:p-norm-on-Schw}
    \|u\|_p := \| \LapS^p u\|_\energy.
  \end{equation}
  Since $(\dom \LapS^p, \|\cdot\|_p)$ is a Hilbert space for each $p \in \bN$, the subspace \Schw is a Fr\'{e}chet space. Note also that since \Lap is unbounded, \Schw is a \emph{proper} subspace of \HE.
\end{remark}

\begin{lemma}\label{thm:vx-in-Schw}
  If $\deg(x)$ is finite for each $x \in \verts$, %or if $\uBd < \iy$, 
  then one has $v_x \in \Schw$.
  \begin{proof}
    %If $\uBd < \iy$, then it is not hard to show that \Lap is bounded, and hence defined everywhere. 
    If $\deg(x) < \iy$ then Lemma~\ref{thm:dx-as-vx} shows that $\gd_x \in \spn\{v_x\}$.
  \end{proof}
\end{lemma}

Note that we take $\deg(x) < \iy$ as a blanket assumption, in fact, as part of the definition of a network. However, it is emphasized in Lemma~\ref{thm:vx-in-Schw} because this is the only place where it is really necessary. (But note that $\deg(x)$ may be unbounded.)

\begin{lemma}\label{thm:dx-as-vx}
  For any $x \in \verts$, $\gd_x = \cond(x) v_x - \sum_{y \nbr x} \cond_{xy} v_y$.
  \begin{proof}
    Lemma~\ref{thm:<delta_x,v>=Lapv(x)} implies $\la \gd_x, u\ra_\energy = \la \cond(x) v_x - \sum_{y \nbr x} \cond_{xy} v_y, u\ra_\energy$ for every $u \in \HE$, so apply this to $u=v_z$, $z \in \verts$. Since $\gd_x, v_x \in \HE$, it must also be that $\sum_{y \nbr x} \cond_{xy} v_y \in \HE$.
  \end{proof}
\end{lemma}

\begin{remark}\label{rem:Lap-maps-span(vx)-into-itself}
  When the hypotheses of Lemma~\ref{thm:vx-in-Schw} are satisfied, it should be noted that $\spn\{v_x\}$ is dense in \Schw with respect to \energy, but \emph{not} with respect to the Fr\'echet topology induced by the seminorms \eqref{eqn:p-norm-on-Schw}, nor with respect to the graph norm. One has the inclusions
  \linenopax
  \begin{equation}\label{eqn:vx-graph-inclusions}
    \left\{\left[\begin{array}{c} v_x \\ \LapM v_x \end{array}\right]\right\}
    \ci \left\{\left[\begin{array}{c} s \\ \LapS s \end{array}\right]\right\}
    \ci \left\{\left[\begin{array}{c} u \\ \LapS u \end{array}\right]\right\}
  \end{equation}
  where $s \in \Schw$ and $u \in \HE$. The second inclusion is dense but the first is not.
\end{remark}

%%!TEX root = bdG.tex

\section{The structure of $\Schw$ and $\Schw'$}
\label{sec:Structure-of-Schw}

From this point on, we assume that a Gel'fand triple has been chosen, using either of the methods described in the previous section. Henceforth, we use the symbol \nuc to denote the operator $\bar\Onb=\Onb^\ad$ or the operator \LapS, depending on how the Gel'fand triple was constructed:
  \linenopax
  \begin{align}\label{eqn:nuc}
    \nuc := 
    \begin{cases}
      \bar\Onb, & \text{Definition~\ref{def:Onb}}\\
      \LapS, & \text{Definition~\ref{def:extn-of-Lap}}.
    \end{cases}
  \end{align}

\subsection{The structure of $\Schw$}
%\label{sec:Structure-of-Schw}

We establish that \Schw is a dense analytic subset of \HE, and that the energy product can be extends not just to a pairing on $\Schw \times \Schw'$, but all the way to a pairing on $\HE \times \Schw'$. Parts of this subsection closely parallel the general theory, and a good reference would be \cite{Simon} or \cite{Hida80}.

\begin{defn}\label{def:spectral-truncation}
  Let $\charfn{[a,b]}$ denote the usual indicator function of the interval $[a,b] \ci \bR$, and let \spectrans be the spectral transform in the spectral representation of \nuc, and let $E$ be the associated projection-valued measure. Then define $E_n$ to be the \emph{spectral truncation operator} acting on \HE by
  \linenopax
  \begin{align*}%\label{eqn:}
     E_n u 
     := \spectrans^\ad \charfn{[\frac1n,n]} \spectrans u 
     = \int_{1/n}^n E(dt)u.
  \end{align*}
\end{defn}

%  Note that Lemma~\ref{thm:vx-in-Schw} implies \Schw is dense in \HE. However, we give an alternative proof in Lemma~\ref{thm:Schw-dense-in-HE}, as this yields formula \eqref{eqn:def:spectral-truncation}, which is useful for approximating elements of \HE by elements of \Schw. In particular, it is used in the proof of Theorem~\ref{thm:HE-isom-to-L2(S',P)}.

\begin{lemma}\label{thm:Schw-dense-in-HE}
  With respect to \energy, \Schw is a dense analytic subspace of \HE.
  \begin{proof}
    This essentially follows immediately once it is clear that $E_n$ maps \HE into \Schw. For $u \in \HE$, and for any $p=1,2,\dots$,
    \linenopax
    \begin{equation}\label{eqn:E-norm-of-Lapp(spectral-truncation)}
      \|\nuc^{p} E_n u\|_\energy^2 
      = \int_{1/n}^n \gl^{2p} \|E(d\gl)u\|_\energy^2
      \leq n^{2p} \|u\|_\energy^2,
    \end{equation}
    So $E_n u \in \Schw$. It follows that $\|u-E_n u\|_\energy \to 0$ by standard spectral theory.
  \end{proof}
\end{lemma}

\begin{theorem}\label{thm:Wiener-product-as-Lap(p)-powers}
  $\Schw \ci \HE \ci \Schw'$ is a Gel'fand triple, and the energy form $\la \cdot,\cdot\ra_\energy$ extends to a pairing on $\Schw \times \Schw'$ defined by
 \linenopax
  \begin{equation}\label{eqn:energy-extends-by-ps}
    \la u,v\ra_\Wiener := \la \nuc^p u, \nuc^{-p} v \ra_\energy,
  \end{equation}
  where $p$ is any integer such that $|v(u)| \leq K \|\Lap^p u\|_\energy$ for all $u \in \Schw$, for some $K>0$.
  \glossary{name={$\la\cdot,\cdot\ra_\Wiener$},description={(Wiener) inner product on $L^2(\Schw',\prob)$},sort=<,format=textbf}
  \begin{proof}
    In combination with \eqref{eqn:Schwartz}--\eqref{eqn:p-norm-on-Schw}, Lemma~\ref{thm:Schw-dense-in-HE} establishes that $\Schw \ci \HE \ci \Schw'$ is a Gel'fand triple.
    %For \eqref{eqn:energy-extends-by-ps}, the definition of $\Schw'$ ensures we can find a $p$ and $K$ to satisfy the estimate; it remains to explain the meaning of $\la \nuc^p u, \nuc^{-p} v \ra_\energy$ for $v \in \Schw'$. 
    If $v \in \Schw'$, then there is a $C$ and $p$ such that $|\la s, v\ra_\Wiener| \leq C\|\nuc^p s\|_\energy$ for all $s \in \Schw$. Set $\gf(\nuc^p s) := \la s, v\ra_\Wiener$ to obtain a continuous linear functional on \HE (after extending to the orthogonal complement of $\spn\{\nuc^p s\}$ by 0 if necessary). Now Riesz's lemma gives a $w \in \HE$ for which $\la s, v\ra_\Wiener = \la \nuc^p s, w\ra_\energy$ for all $s \in \Schw$ and we define $\nuc^{-p} v := w \in \HE$ to make the meaning of the right-hand side of \eqref{eqn:energy-extends-by-ps} clear. 
  \end{proof}
\end{theorem}

\begin{lemma}\label{thm:energy-extends}
 The pairing on $\Schw \times \Schw'$ is equivalently given by
  \linenopax
  \begin{equation}\label{eqn:energy-extends-by-vn}
    \la u,\gx\ra_\Wiener = \lim_{n \to \iy} \gx(E_n u), %\la u, \tilde E_n \gx \ra_\energy,
  \end{equation}
  where the limit is taken in the topology of $\Schw'$. Moreover, $\tilde u(\gx) = \la u, \gx\ra_\Wiener$ is \bR-valued on $\Schw'$.
  \begin{proof}
     %We claim that $\nuc^{-p} \tilde E_n = E_n \nuc^{-p}$. It suffices to check this on vectors of the form $\nuc^p s$:
    %\linenopax
    %\begin{align*}
    %  \la \nuc^p s, \nuc^{-p} \tilde E_n\gx \ra_\energy 
    %  = \la s, \tilde E_n \gx \ra_\energy 
    %  = \la E_n s, \gx \ra_\energy
    %  \qq s \in \Schw,
    %\end{align*}
    %by the definition of $\tilde E_n$ from Lemma~\ref{thm:Spectruncation-extends-to-S'}, and 
    $E_n$ commutes with \nuc. This is a standard result in spectral theory, as $E_n$ and \nuc are unitarily equivalent to the two commuting operations of truncation and multiplication, respectively. Therefore,
    \linenopax
    \begin{align*}
      \gx(E_n u)
      = \la E_n u, \gx\ra_\Wiener
      = \la \nuc^p E_n s, \nuc^{-p} \gx \ra_\energy
      = \la E_n \nuc^p s, \nuc^{-p} \gx \ra_\energy
      = \la\nuc^p s, E_n \nuc^{-p} \gx \ra_\energy.
    \end{align*}
    Standard spectral theory also gives $E_n v \to v$ in \HE, so 
    \linenopax
    \begin{align*}
      \lim_{n \to \iy} \gx(E_n u)
      = \lim_{n \to \iy} \la\nuc^p s, E_n \nuc^{-p} \gx \ra_\energy
      = \la \nuc^p u,\nuc^{-p} v \ra_\energy.
    \end{align*}
    %where the limits are taken in the topology of $\Schw'$.
    %To check that \eqref{eqn:energy-extends-by-vn} is well-defined, use the definition of $\Schw'$ to find a $C$ and $p$ such that
    %\begin{align*}%\label{eqn:}
    %  |\gx(E_n u)|
    %  \leq C \|\nuc^p E_n u\|_\energy
    %  \leq C n^p \|u\|_\energy,
    %\end{align*}
    %where the latter bound comes by \eqref{eqn:E-norm-of-Lapp(spectral-truncation)} again. Then Riesz's lemma gives $w_n \in \HE$ for which $\gx(E_n u) = \la v, w_n\ra_\energy$.
    
    Note that the pairing $\la \cdot \,,\,\cdot \ra_\Wiener$ is a limit of real numbers, and hence is real.
  \end{proof}
\end{lemma}
%   then $\Lap^{-p} v$ may not lie in \HE. However, $\tilde E_n v \in \HE$ by Lemma~\ref{thm:Spectruncation-extends-to-S'}, whence $\Lap^{-p} \tilde E_n v \in \HE$ because the spectral resolution implies $\|\Lap^{-p} E_n v\|_\energy^2 \leq n^{2p} \|v\|_\energy^2$ by the same computation as in \eqref{eqn:E-norm-of-Lapp(spectral-truncation)}.
    
    %Since \nuc is self-adjoint, we understand $\la u, v_n\ra_\Wiener = \la \nuc^p u, \nuc^{-p} v_n \ra_\energy$. 
    
\begin{cor}\label{thm:Spectruncation-extends-to-S'}
  $E_n$ extends to a mapping $\tilde E_n: \Schw' \to \HE$ defined via $\la u, \tilde E_n \gx \ra_\energy := \gx(E_n u)$. Thus, we have a pointwise extension of $\la \cdot \,,\,\cdot\ra_\Wiener$ to $\HE \times \Schw'$ given by
  \begin{equation}\label{eqn:energy-extends-by-truncating-on-S'}
    \la u,\gx\ra_\Wiener = \lim_{n \to \iy} \la u, \tilde E_n \gx \ra_\energy.
  \end{equation}
\end{cor}

%    For a \bR-valued $v \in \HE$, one can find $\{v_n\} \ci \Schw$ with $\|v-v_n\|_\energy \to 0$ by choosing the spectral truncations $v_n := E_n v$. Since $\cj{\nuc u} = \nuc \cj{u}$, the spectral theorem gives
%  \begin{align*}%\label{eqn:}
%     \cj{v_n} 
%     = \cj{\int_{1/n}^n E(dt)v} 
%     = \int_{1/n}^n E(dt)\cj{v}
%     = \int_{1/n}^n E(dt)v 
%     = v_n.
%  \end{align*}
%  Since $v \in \Schw'$ is also \bR-valued (cf.~Remark~\ref{rem:Schwartz-space-is-real-valued}), $\la u,v\ra_\Wiener = \lim_{n \to \iy} \la u_n,v\ra_\Wiener$.   

%  \version{}{\marginpar{This is no longer true, since that theorem only has one implication, now.}It is also clear that the dipoles $v_x$ are \emph{not} analytic vectors. If they were, then because they are dense in \HE, a theorem of Nelson would imply that \LapM has a unique self-adjoint extension. However, we know this to be false in general, by Theorem~\ref{thm:LapV-not-ess-selfadjoint-iff-Harm=0}. Consequently, one is led to expect the poles of $v_x$ to correspond to the boundary $\bd \Graph$ somehow; the authors are currently investigating this matter.}

\begin{comment}\label{thm:vx-in-Schw}
  If every vertex of \Graph has finite degree, then
  $\Fin \ci V=\spn\{v_x\}_{x \in \verts}$.
  \begin{proof}
    Let $x \neq o$ and note that Lemma~\ref{thm:<delta_x,v>=Lapv(x)} gives
    \linenopax
    \begin{align*}
      \gd_z(x) 
      = \gd_x(z) - \gd_x(o)
      = \la \gd_x, v_z\ra_\energy
      &= \Lap v_z(x) = \cond(x) v_z(x) - \sum_{y \nbr x} \cond_{xy} v_z(y).
      &&\qedhere
    \end{align*}
  \end{proof}
\end{comment}

\subsection{The structure of $\Schw'$}
\label{sec:Structure-of-Schw'}

The next results are structure theorems akin to those found in the classical theory of distributions; see \cite[\S6.3]{Str03} or \cite[\S3.5]{Al-Gwaiz}. If $\HE \ci \Schw$, then Theorem~\ref{thm:structure-of-S'} would say $\Schw' = \bigcup_p \nuc^p(\HE)$ (of course, this is typically false when $\Harm \neq 0$).

\begin{theorem}\label{thm:structure-of-S'}
  The distribution space $\Schw'$ is
  \linenopax
  \begin{align}\label{eqn:structure-of-S'}
    \Schw' = \{\gx(u) = \la \nuc^p u, v\ra_\energy \suth \exists v \in \HE, p \in \bZ^+, \forall u \in \Schw\}.
   \end{align}
  \begin{proof}
    It is clear from the Schwarz inequality that $\gx(u) = \la \nuc^p u, v\ra_\energy$ defines a continuous linear functional on \Schw, for any $v \in \HE$ and nonnegative integer $p$. For the other direction, we use the same technique as in Lemma~\ref{thm:Wiener-product-as-Lap(p)-powers}. Observe that if $\gx \in \Schw'$, then there exists $K,p$ such that $|\gx(u)| \leq K \|\nuc^p u\|_\energy$ for every $u \in \Schw$. This implies that the map $\gx: \nuc^p u \mapsto \gx(u)$ is continuous on the subspace $Y = \spn\{\nuc^p u \suth u \in \HE, p \in \bZ^+\}$. This can be extended to all of \HE by precomposing with the orthogonal projection to $Y$. Now Riesz's lemma gives a $v \in \HE$ for which $\gx(u) = \la \nuc^p u, v\ra_\energy$. % = \la u, \nuc^p v\ra_\Wiener = \la u, \nuc^p f\ra_\Wiener $, where $f = \Pfin v$, whence $\gx = \nuc^p f$.
%    
%    ($\ce$) Let $v \in \Fin$ and fix any $p \in \bZ^+$. Fix $u \in \Schw$ and define $\gx(u) := \la u, \nuc^p v\ra_\Wiener$. Then the bound 
%  \linenopax
%  \begin{align}\label{eqn:bddness-of-F(u)}
%    | \gx(u)| 
%    = \left|\la u, \nuc^p v\ra_\Wiener\right|
%    = \left|\la \nuc^p u, v\ra_\energy\right|
%    \leq \| \nuc^p u\|_\energy \| v \|_\energy,
%  \end{align}
%  shows \gx is a well-defined element of $\Schw'$ corresponding to $\nuc^p v$.
  \end{proof}
\end{theorem}

Note that $v \in \HE$ may not lie in the domain of $\nuc^p$. If it did, one would have $\la \nuc^p u, v\ra_\energy = \la u, \nuc^p v\ra_\Wiener = \la u, \nuc^p f\ra_\Wiener$, where $f = \Pfin v$. The theorem could then be written $\Schw' = \bigcup_{p=0}^\iy \nuc^p (\Fin)$. However, this turns out to have contradictory implications.

We now provide two results enabling one to recognize certain elements of $\Schw'$.

\begin{lemma}\label{thm:Schw-dist-criterion}
  A linear functional $f:\Schw \to \bC$ is an element of $\Schw'$ if and only if there exists $p \in \bN$ and $F_0,F_1,\dots F_p \in \HE$ such that
  \linenopax
  \begin{equation}\label{eqn:Schw-dist-criterion}
    f(u) = \sum_{k=0}^p \la F_k, \nuc^k u \ra_\energy, \q\forall u \in \HE.
  \end{equation}
  \begin{proof}
    By definition, $f \in \Schw'$ iff $\exists p,C<\iy$ for which $|f(u)| \leq C \|u\|_p$ for every $u \in \Schw$. Therefore, the linear functional
    \linenopax
    \begin{align*}
      \gF:\bigoplus\nolimits_{k=0}^p \dom(\nuc^k) \to \bC
      \qq\text{by}\qq
      \gF(u,\nuc u, \nuc^2 u, \dots \nuc^p u) = f(u)
    \end{align*}
    is continuous and Riesz's Lemma gives $F = (F_k)_{k=0}^p \in \bigoplus_{k=0}^p \negsp[5]\HE$ with
    \linenopax
    \begin{align*}
      f(u) &= \la F, (u,\nuc u, \dots \nuc^p u)\ra_{\bigoplus \HE}
        = \sum_{k=0}^p \la F_k, \nuc^k u \ra_{\bigoplus\negsp[5]\HE}.
        &&\qedhere
    \end{align*}
  \end{proof}
\end{lemma}

\begin{cor}\label{thm:Schw'=HE-when-Lap-bdd}
  If $\nuc:\HE \to \HE$ is bounded, then $\Schw'=\HE$.
  \begin{proof}
    We always have the inclusion $\HE \hookrightarrow \Schw'$ by taking $p=0$. If  \nuc is bounded, then the adjoint $\nuc^\ad$ is also bounded, and \eqref{eqn:Schw-dist-criterion} gives
    \linenopax
    \begin{equation}\label{eqn:Schw'=HE-when-Lap-bdd}
      f(u) = \left\la \sum_{k=0}^p (\nuc^\ad)^k  F_k, u \right\ra_{\bigoplus\negsp[5]\HE},
      \q\forall u \in \Schw.
    \end{equation}
    Since \Schw is dense in \HE by Lemma~\ref{thm:Schw-dense-in-HE}, we have $f = \sum_{k=0}^p (\nuc^\ad)^k  F_k \in \HE$.
  \end{proof}
\end{cor}

\begin{remark}\label{rem:unbounded-iff-infinite}
  In view of Lemma~\ref{thm:properties-of-Onb}, Corollary~\ref{thm:Schw'=HE-when-Lap-bdd} shows that $\Schw'$ is a proper extension of \HE on any infinite network.
\end{remark}

In the case when the Gel'fand triple is constructed from the domain of \LapM, as in Definition~\ref{def:extn-of-Lap}, then one can extend \Lap to distributions.

\begin{defn}\label{def:Lap-on-distributions}
  Extend \Lap to $\Schw'$ by defining
  \linenopax
  \begin{equation}\label{eqn:Lap-on-distributions}
    \Lap \gx(v_x) := \la \gd_x, \gx \ra_\Wiener,
  \end{equation}
  so that $\Lap\gx(v_x) = \sum_{y \nbr x} \cond_{xy} (\gx(v_x)-\gx(v_y))$ follows readily from Lemma~\ref{thm:dx-as-vx}. 
  
  Now extend \Lap to $\tilde \Lap$ defined on $\tilde v_x \in L^2(\Squoth,\quotprob)$ by $\tilde \Lap (\tilde v_x)(\gx) := \widetilde{\Lap v_x} (\gx)$, so that 
  \linenopax
  \begin{equation}\label{eqn:Lap-on-L2(distns)}
    \tilde \Lap: \tilde v_x \mapsto \cond(x) \tilde v_x - \sum_{y \nbr x} \cond_{xy} \tilde v_y.
  \end{equation}
  Since $v_x \mapsto \tilde v_x$ is an isometry, it is no great surprise that
  \linenopax
  \begin{align}\label{eqn:Lap-on-distributions-inner-product}
    \la \tilde v_x, \tilde \Lap \tilde v_y\ra_{L^2}
    &= \int_{\Schw'} \tilde v_x(\gx) \tilde v_y(\Lap \gx)\,d\quotprob(\gx)
    = \la v_x, \Lap v_y\ra_\energy.
  \end{align}
\end{defn}

%%!TEX root = bdG.tex

\section{The Wiener embedding and the space $\Schw'$}
\label{sec:Wiener-embedding}

We have now obtained a Gel'fand triple $\Schw \ci \HE \ci \Schw'$ (from either Lemma~\ref{thm:sandwich=triple} or Theorem~\ref{thm:Wiener-product-as-Lap(p)-powers}), and we are ready to apply the Minlos Theorem to a particularly lovely positive definite function on \Schw, in order that we may obtain a particularly nice measure on $\Schw'$. This allows us to realize $\bd \Graph$ as a subset of $\Schw'$. Recall that \Schw contains the energy kernel; see Lemma~\ref{thm:v_x-in-S} or Lemma~\ref{thm:vx-in-Schw}.

\subsection{The Wiener embedding}

In \cite[\S5]{ERM}, we constructed \HE from the resistance metric by making use of negative definite functions. We now apply this to a famous result of Schoenberg which may be found in \cite{Ber84,ScWh49}.
\begin{theorem}[Schoenberg]
  \label{thm:Schoenberg's-Thm}
  Let $X$ be a set and let $Q: X \times X \to \bR$ be a function. Then the following are equivalent.
  \begin{enumerate}
    \item $Q$ is negative definite.
    \item $\forall t \in \bR^+$, the function $p_t(x,y) := e^{-t Q(x,y)}$ is positive definite on $X \times X$.
    \item There exists a Hilbert space $\sH$ and a function $f:X \to \sH$ such that $Q(x,y) = \|f(x)-f(y)\|_\sH^2$.
  \end{enumerate}
\end{theorem}
       \begin{comment}
         Euclidean space $\bRd$ is self-dual. Therefore, the Fourier transform $Q=Q_\gn$ of any given positive finite Borel measure \gn on \bRd is the continuous function
         \linenopax
         \begin{align*}
           Q(t) = \int_{\bRd} e^{-\ii t x} \,d\gn(x).
         \end{align*}
         $Q$ is continuous because $e^{-\ii t x}$ is a continuous and periodic function of $t$ for any fixed $x$. Observe also that the kernel $K(x,y) = Q(y-x)$ is positive definite, whence $Q$ is a positive definite function.
       \end{comment}

In the proof of the following theorem, we apply Schoenberg's Theorem with $t=\frac12$ to the resistance metric in the form 
\linenopax
\begin{align}\label{eqn:def:R^F(x,y)-energy}
  R^F(x,y) = \|v_x-v_y\|_\energy^2,
\end{align}
which appears in \cite[Thm.~2.13]{ERM}.
The proof of Theorem~\ref{thm:HE-isom-to-L2(S',P)} also uses the notation $\Ex_\gx(f) := \int_{\Schw'} f(\gx) \,d\prob(\gx)$. 
\glossary{name={$\Ex_\gx(f)$},description={expectation of $f$ with dummy variable \gx; $\int f(\gx)\,d\gm(x)$},sort=E,format=textbf}

\begin{theorem}[Wiener embedding]
  \label{thm:HE-isom-to-L2(S',P)}
  The Wiener transform $\sW:\HE \to L^2(\Schw',\prob)$ by
  \linenopax
  \begin{equation}\label{eqn:Gaussian-transform}
    \sW : v \mapsto \tilde v,
    \q \tilde v(\gx) := \la v, \gx\ra_\Wiener,
  \end{equation}
  is an isometry. The extended reproducing kernel $\{\tilde v_x\}_{x \in \verts}$ is a system of Gaussian random variables which gives the resistance distance by
  \linenopax
  \begin{equation}\label{eqn:R(x,y)-as-expectation}
    R^F(x,y) = \Ex_\gx((\tilde v_x - \tilde v_y)^2).
  \end{equation}
  Moreover, for any $u,v \in \HE$, the energy inner product extends directly as
  \linenopax
  \begin{equation}\label{eqn:Expectation-formula-for-energy-prod}
    \la u, v \ra_\energy
    = \Ex_\gx\left( \cj{\tilde{u}} \tilde{v} \right)
    = \int_{\Schw'} \cj{\tilde{u}} \tilde{v} \,d\prob.
  \end{equation}
  \begin{proof}
    Since $R^F(x,y)$ is negative semidefinite (see \cite[Thm.~5.4]{ERM}), we may apply Schoenberg's theorem and deduce that $\exp(-\tfrac12\|u-v\|_\energy^2)$ is a positive definite function on $\HE \times \HE$. Consequently, an application of the Minlos correspondence to the Gel'fand triple established in Lemma~\ref{thm:Schw-dense-in-HE} yields a Gaussian probability measure \prob on $\Schw'$. 
    
    Moreover, \eqref{eqn:Minlos-identity} gives 
    \linenopax
    \begin{align}\label{eqn:Minlos-eqns}
      \Ex_\gx(e^{\ii\la u, \gx \ra_\Wiener}) = e^{-\frac12\|u\|_\energy^2},
    \end{align}
    whence one computes
    \linenopax
    \begin{align}\label{eqn:Minlos-expectation-integral}
      \int_{\Schw'} \left(1 + \ii\la u, \gx\ra_\Wiener - \frac12\la u, \gx\ra_\Wiener^2 + \cdots \right)\,d\prob(\gx) 
      = 1 - \frac12 \la u, u \ra_\energy + \cdots.
    \end{align}
    Now it follows that $\Ex(\tilde{u}^2) = \Ex_\gx(\la u, \gx\ra_\Wiener^2) = \|u\|_\energy^2$ for every $u \in \Schw$, by comparing the terms of \eqref{eqn:Minlos-expectation-integral} which are quadratic in $u$. Therefore, $\sW:\HE \to \Schw'$ is an isometry, and \eqref{eqn:Minlos-expectation-integral} gives
    \linenopax
    \begin{align}\label{eqn:Exp-tilde-vx=E(vx)}
      \Ex_\gx(|\tilde v_x - \tilde v_y|^2)
      = \Ex_\gx(\la v_x - v_y, \gx \ra^2)
      &= \|v_x - v_y\|_\energy^2,
    \end{align}
    whence \eqref{eqn:R(x,y)-as-expectation} follows from \eqref{eqn:def:R^F(x,y)-energy}. Note that by comparing the linear terms, \eqref{eqn:Minlos-expectation-integral} implies $\Ex_\gx(1) = 1$, so that \prob is a probability measure, and $\Ex_\gx(\la u,\gx\ra) = 0$ and $\Ex_\gx(\la u,\gx\ra^2) = \|u\|_\Wiener^2$, so that \prob is actually Gaussian.
    
    %Approximating $v_x$ by a Cauchy sequence $\{v_x^{(n)}\}_{n=1}^\iy \ci \Schw$ as in Lemma~\ref{thm:Schw-dense-in-HE}, we see that $\tilde v_x$ as in \eqref{eqn:Gaussian-transform} is a well defined element of $L^2(\Schw',\prob)$. Fisher's Theorem gives a subsequence of $\{v_x^{(n)}\}$ which converges pointwise \prob-a.e. to $\tilde v_x$, and so $\tilde v_x$ is a bona fide random variable.
    
    Finally, use polarization to compute
    \linenopax
    \begin{align*}
      \la u, v \ra_\energy
      &= \frac14 \left(\|u+v\|_\energy^2 - \|u-v\|_\energy^2\right) \\
      &= \frac14 \left(\Ex_\gx\left(\left|\tilde{u}+ \tilde{v}\right|^2\right) 
        - \Ex_\gx\left(\left|\tilde{u}-\tilde{v}\right|^2\right)\right) 
        &&\text{by \eqref{eqn:Exp-tilde-vx=E(vx)}} \\
      &= \frac14 \int_{\Schw'} \left|\tilde{u}+\tilde{v}\right|^2(\gx)
        - \left|\tilde{u}-\tilde{v}\right|^2(\gx) \,d\prob(\gx) \\
      &= \int_{\Schw'} \cj{\tilde{u}}(\gx) \tilde{v}(\gx) \,d\prob(\gx).
    \end{align*}
    This establishes \eqref{eqn:Expectation-formula-for-energy-prod} and completes the proof.
  \end{proof}
\end{theorem}

  It is important to note that since the Wiener transform $\sW : \Schw \to \Schw'$ is an isometry, the conclusion of Minlos' theorem is stronger than usual: the isometry allows the energy inner product to be extended isometrically to a pairing on $\HE \times \Schw'$ instead of just $\Schw \times \Schw'$.
  
\begin{comment}\label{thm:pairing-extends-to-HE}
  The pairing $\la \cdot , \cdot \ra_\Wiener$ is well-defined on $\HE \times \Schw'$.
  \begin{proof}
    For $u \in \HE$, let $\{u_n\} \ci \Schw$ be a sequence such that $\|u_n-u\|_\energy \to 0$. Then $\{u_n\}$ is Cauchy in \HE and then Theorem~\ref{thm:HE-isom-to-L2(S',P)} gives
    \linenopax
    \begin{align*}
      \|\widetilde{u_n} - \widetilde{u_{n+m}}\|_\Wiener 
      = \|{u_n} -{u_{n+m}}\|_\energy,
    \end{align*}
    so $\{\widetilde{u_n}\}$ is Cauchy in $L^2(\Schw',\prob)$, and hence has a limit $\tilde u$. Then Fischer's theorem gives $\tilde u(\gx) = \lim_{n \to \iy} \widetilde{u_{n}}(\gx)$ for \prob-a.e. $\gx \in \Schw'$ (passing to a subsequence if necessary), whence
    \linenopax
    \begin{align*}
      \la u, \gx \ra_\Wiener
      := \lim_{n \to \iy} \la u_n, \gx \ra_\Wiener
    \end{align*}
    gives a well-defined extension of the pairing to $\HE \times \Schw'$.
  \end{proof}
\end{comment}

\begin{remark}\label{rem:HE-into-L2(S',prob)-gives-hermitian-multiplication}
  With the embedding $\HE \to L^2(\Schw',\prob)$, we obtain a maximal abelian algebra of Hermitian multiplication operators $L^\iy(\Schw')$ acting on $L^2(\Schw',\prob)$. For a sharp contrast, note that the multiplication operators on \HE are trivial, by \cite[Lem.~2.3]{DGG}. This result states that if $\gf:\verts \to \bR$ and $M_\gf$ denotes the multiplication operator defined by $(M_\gf u)(x) = \gf(x) u(x)$, then $M_\gf$ is Hermitian if and only if $M_\gf = k\id$, for some $k \in \bR$.
\end{remark}

\begin{remark}\label{rem:complexify-L2(S',P)}
  The reader will note that we have taken pains to keep everything \bR-valued in this section (especially the elements of \Schw and $\Schw'$), primarily to ensure the convergence of $\int_{\sS'} e^{\ii\la u, \gx \ra_\Wiener}\,d\prob(\gx)$ in \eqref{eqn:Minlos-eqns}. However, now that we have established the fundamental identity $\la u, v \ra_\energy = \int_{\sS'} \cj{\tilde{u}} \tilde{v} \,d\prob$ in
\eqref{eqn:Expectation-formula-for-energy-prod} and extended the pairing $\la \cdot , \cdot \ra_\Wiener$ to $\HE \times \Schw'$, we are at liberty to complexify our results via the standard decomposition into real and complex parts: $u = u_1 + \ii u_2$ with $u_i$ \bR-valued elements of \HE, etc.
\end{remark}

\begin{remark}\label{rem:Wiener-improves-Minlos}
  The polynomials are dense in $L^2(\Schw',\prob)$. More precisely, if $\gf(t_1, t_2, \dots, t_k)$ is an ordinary polynomial in $k$ variables, then 
  \linenopax
  \begin{align}\label{eqn:polynomials-in-S'}
  \gf(\gx) := \gf\left(\la u_1, \gx\ra_\Wiener, \la u_2, \gx\ra_\Wiener, \dots \vstr[2.2] \, \la u_k,\gx\ra_\Wiener\right)
  \end{align}
  is a polynomial on $\Schw'$ and %(with complex coefficients), then 
  \linenopax
  \begin{align}\label{eqn:Poly(n)-in-S'}
    \Poly_n := \{\gf\left(\widetilde{u_1}(\gx), \widetilde{u_2}(\gx), \dots \vstr[2.2] \, \widetilde{u_k}(\gx)\right), \deg(\gf) \leq n, \suth u_j \in \HE, \gx \in \Schw'\}
  \end{align}
  is the collection of polynomials of degree at most $n$, and $\{\Poly_n\}_{n=0}^\iy$ is an increasing family whose union is all of $\Schw'$. One can see that the monomials $\la u, \gx\ra_\Wiener$ are in $L^2(\Schw',\prob)$ as follows: compare like powers of $u$ from either side of \eqref{eqn:Minlos-expectation-integral} to see that $\Ex_\gx\left(\la u, \gx\ra_\Wiener^{2n+1}\right) = 0$ and
  \linenopax
  \begin{align}\label{eqn:expectation-of-even-monomials}
    \Ex_\gx\left(\la u, \gx\ra_\Wiener^{2n}\right) 
    = \int_{\Schw'} |\la u, \gx\ra_\Wiener|^{2n} \, d\prob(\gx) 
    = \frac{(2n)!}{2^n n!} \|u\|_\energy^{2n}, 
  \end{align}
  and then apply the Schwarz inequality. 
  
  To see why the polynomials $\{\Poly_n\}_{n=0}^\iy$ should be dense in $L^2(\Schw',\prob)$ observe that the sequence $\{P_{\Poly_n}\}_{n=0}^\iy$ of orthogonal projections increases to the identity, and therefore, $\{P_{\Poly_n} \tilde u\}$ forms a martingale, for any $u \in \HE$ (i.e., for any $\tilde u \in L^2(\Schw',\prob)$).
  
  Denote the ``multiple Wiener integral of degree $n$'' by 
  \linenopax
  \begin{align*}
    H_n := %\Poly_n - \Poly_{n-1} = 
    \left(cl \spn\{\la u, \cdot\ra_\Wiener^n \suth u \in \HE\}\right) \ominus \{\la u, \cdot\ra_\Wiener^k \suth k<n, u \in \HE\},
  \end{align*}
  for each $n \geq 1$, and $H_0 := \bC \one$ for a vector \one with $\|\one\|_2=1$. Then we have an orthogonal decomposition of the Hilbert space     
  \linenopax
  \begin{align}\label{eqn:Fock-repn-of-L2(S',P)}
    L^2(\Schw',\prob) = \bigoplus_{n=0}^\iy H_n.
  \end{align}
  See \cite[Thm.~4.1]{Hida80} for a more extensive discussion.
  
  \glossary{name={\one},description={the constant function 1, the vacuum vector},sort=1,format=textbf}
A physicist would call \eqref{eqn:Fock-repn-of-L2(S',P)} the Fock space representation of $L^2(\Schw',\prob)$ with ``vacuum vector'' \one; note that $H_n$ has a natural (symmetric) tensor product structure. Familiarity with these ideas is not necessary for the sequel, but the decomposition \eqref{eqn:Fock-repn-of-L2(S',P)} is helpful for understanding two key things:
  \begin{enumerate}[(i)]
    \item The Wiener isometry $\sW:\HE \to L^2(\Schw',\prob)$ identifies \HE with the subspace $H_1$ of $L^2(\Schw',\prob)$, in particular, $L^2(\Schw',\prob)$ is not isomorphic to \HE. In fact, it is the second quantization of \HE.
    \item The constant function \one is an element of $L^2(\Schw',\prob)$ but does not correspond to any element of \HE. In particular, constant functions in \HE are equivalent to 0, but this is not true in $L^2(\Schw',\prob)$.
  \end{enumerate}
  It is somewhat ironic that we began this story by removing the constants (via the introduction of \energy), only to reintroduce them with a certain amount of effort, much later. Item (ii) explains why it is not nonsense to write things like $\prob(\Schw') = \int_{\Schw'} \one \,d\prob = 1$, and will be helpful when discussing boundary elements in \S\ref{sec:bdG-as-equivalence-classes-of-paths}.
\end{remark}

\begin{cor}\label{thm:resistance-as-distributional-integral}
  For $e_x(\gx) := e^{\ii\la v_x, \gx\ra_\Wiener}$, one has $\Ex_\gx(e_x) = e^{-\frac12 R^F(o,x)}$ and hence
  \linenopax
  \begin{equation}\label{eqn:resistance-as-distributional-integral}
    \Ex_\gx(\cj{e_x} e_y) 
    = \int_{\Schw'} \cj{e_x(\gx)} e_y(\gx)\,d\prob
    = e^{-\frac12 R^F(x,y)}.
  \end{equation}
  \begin{proof}
    Substitute $u=v_x$ or $u=v_x-v_y$ in \eqref{eqn:Minlos-eqns} and apply \eqref{eqn:def:R^F(x,y)-energy}.
  \end{proof}
\end{cor}

\begin{remark}
  Free resistance is interpreted as the reciprocal of an integral over a path space in \cite[Rem.~3.15]{ERM}; Corollary~\ref{thm:resistance-as-distributional-integral} provides a variation on this theme:
  \linenopax
  \begin{equation}\label{eqn:resistance-as-distributional-integral}
    R^F(x,y) = -2 \log \Ex_\gx(\cj{e_x} e_y) = -2 \log \int_{\Schw'} \cj{e_x(\gx)} e_y(\gx)\,d\prob.
  \end{equation}
  
  Observe that Theorem~\ref{thm:HE-isom-to-L2(S',P)} was carried out for the free resistance, but all the arguments go through equally well for the wired resistance; note that $R^W$ is similarly negative semidefinite by Theorem~\ref{thm:Schoenberg's-Thm} and \cite[Cor.~5.5]{ERM}. Thus, there is a corresponding Wiener transform $\sW:\Fin \to L^2(\Schw',\prob)$ defined by
  \linenopax
  \begin{equation}\label{eqn:wired-Gaussian-transform}
    \sW : v \mapsto \tilde f,
    \qq f = \Pfin v \;\text{ and }\; \tilde f(\gx) = \la f, \gx\ra_\Wiener.
  \end{equation}
  Again, $\{\tilde f_x\}_{x \in \verts}$ is a system of Gaussian random variables which gives the wired resistance distance by $R^W(x,y) = \Ex_\gx((\tilde f_x - \tilde f_y)^2)$. %The authors are presently working to see if there is any relation between these measures and the free and wired uniform spanning forest measures of \cite[\S9]{Lyons:ProbOnTrees}.
\end{remark}

\begin{remark}
  \label{rem:abuse-of-extension-notation}
  For $u \in \Harm$ and $\gx \in \Schw'$, let us abuse notation and write $u$ for $\tilde{u}$. That is, $u(\gx) := \tilde{u}(\gx) = \la u, \gx\ra_\Wiener$. Unnecessary tildes obscure the presentation and the similarities to the Poisson kernel in \S\ref{sec:Operator-theoretic-interpretation-of-bdG}.
\end{remark}

\begin{remark}\label{rem:onb-gives-Gaussian-iid}
  Theorem~\ref{thm:HE-isom-to-L2(S',P)} showed that $\{\tilde \onb_x\}$ forms a system of Gaussian random variables.
  Since the Wiener transform is an isometry, 
  \linenopax
  \begin{align}\label{eqn:onb-moments}
    \Ex(\tilde \onb_x) = 0  
    \qq\text{and}\qq
    \Ex(\tilde \onb_x \tilde \onb_y) = \gd_{x,y}.
  \end{align}
  Since independence of \emph{Gaussian} random variables is determined by the first two moments, it follows that $\{\tilde \onb_x\}$ forms a system of i.i.d. Gaussian random variables with mean $0$ and variance $1$. This is noteworthy because while independence implies orthogonality, the converse does not hold without the additional hypothesis that the distributions be Gaussian. 
\end{remark}

%%!TEX root = bdG.tex

\section{The resistance boundary of a transient network}
\label{sec:Operator-theoretic-interpretation-of-bdG}

With the tools developed in \S\ref{sec:Gel'fand-triple-for-HE} and \S\ref{sec:Wiener-embedding}, we now construct the resistance boundary $\bd \Graph$ as a set of equivalence classes of infinite paths. Recall that we began with a comparison of the Poisson boundary representation for bounded harmonic functions with the boundary sum representation recalled in \eqref{eqn:bdy-repn-preview}:
\linenopax
\begin{equation*}%\label{eqn:Poisson-bdy-repn-recall}
  u(x) = \int_{\del \gW} u(y) k(x,dy) 
  \qq\leftrightarrow\qq
  u(x) = \sum_{\bd \Graph} u \dn{h_x} + u(o).
\end{equation*}
In this section, we replace the sum with an integral and complete the parallel.

\begin{comment}
\begin{defn}\label{def:quotprob}
  Define $\quotprob$ to be the image measure on $\Schw'/\Fin$ induced by the standard projection $\gp: \Schw' \to \Schw'/\Fin$, i.e., $\quotprob(B) := \prob(\gp^{-1} B)$, for $B \in \sB(\Schw'/\Fin)$. 
\end{defn}
  \glossary{name={$\quotprob$},description={(quotient) probability measure induced by canonical projection},sort=P,format=textbf}

Now \quotprob is a probability measure on the quotient and Theorem~\ref{thm:HE-isom-to-L2(S',P)} gives a corresponding energy integral representation.
\end{comment}

\begin{cor}[Boundary integral representation for harmonic functions]
  \label{thm:Boundary-integral-repn-for-harm} \hfill \\
  For any $u \in \Harm$ and with $h_x = \Phar v_x$, 
  \linenopax
  \begin{equation}\label{eqn:integral-boundary-repn-of-h}
    u(x) = \int_{\Squoth} u(\gx) h_x(\gx) \, d\quotprob(\gx) + u(o).
  \end{equation}
  \begin{proof}
    Starting with \eqref{eqn:def:vx}, compute
    \linenopax
    \begin{align}\label{eqn:boundary-repn-for-harmonic-integral-derived}
      u(x) - u(o)
      = \la h_x, u\ra_\energy
      = \cj{\la u, h_x \ra_\energy}
      = \cj{\int_{\Schw'} \cj{u} h_x \, d\quotprob},
    \end{align}
    where the last equality comes by substituting $v=h_x$ in \eqref{eqn:Expectation-formula-for-energy-prod}. It is shown in \cite[Lem.~2.24]{DGG} that $\cj{h_x}=h_x$. %Note that we are suppressing tildes as in Remark~\ref{rem:abuse-of-extension-notation}. Recall from the Gel'fand triple construction that $\Fin \ci \HE \ci \Schw'$, but note that $h_x(f) = \la h_x, f \ra_\energy = 0 $ for every $f \in \Fin$, so the domain of integration passes to the quotient $\Squoth $. 
  \end{proof}
\end{cor}

\begin{remark}[A Hilbert space interpretation of {bd}\,\Graph]
  \label{rem:boundary-of-G-from-Minlos}
  In view of Corollary~\ref{thm:Boundary-integral-repn-for-harm}, we are now able to ``catch'' the boundary between \Schw and $\Schw'$ by using \nuc and its adjoint. The boundary of \Graph may be thought of as (a possibly proper subset of) $\Squoth$. 
%  \linenopax
%  \begin{equation}\label{eqn:def:boundary-of-G-from-Minlos}
%    \bd \Graph = \Schw' / \Fin.
%  \end{equation}
  Corollary~\ref{thm:Boundary-integral-repn-for-harm} suggests that $\mathbbm{k}(x,d\gx) := h_x(\gx) d\quotprob$ is the discrete analogue in \HE of the Poisson kernel $k(x,dy)$, and comparison of \eqref{eqn:bdy-repn-preview} with \eqref{eqn:integral-boundary-repn-of-h} gives a way of understanding a boundary integral as a limit of Riemann sums:
  \linenopax
  \begin{equation}\label{eqn:boundary-of-G-from-Minlos}
    \int_{\Schw'} u \, h_x \, d\quotprob
    = \lim_{k \to \iy} \sum_{\bd G_k} u(x) \dn{h_x}(x). 
  \end{equation}
  (We continue to omit the tildes as in Remark~\ref{rem:abuse-of-extension-notation}.) By a theorem of Nelson, \quotprob is fully supported on those functions which are H\"{o}lder-continuous with exponent $\ga=\frac12$, which we denote by $\Lip(\frac12) \ci \Schw'$; see \cite{Nelson64}. Recall from \cite[Cor.~2.16]{ERM} that $\HE \ci Lip(\frac12)$.
  See \cite{Arv76a,Arv76b,Minlos63,Nel69}.
\end{remark}

%\input{paths}
%%!TEX root = bdG.tex

\section{Examples}
\label{sec:examples}

Our presentation of $\bd \Graph$ may appear somewhat abstract in the general case. However, we now illustrate the concept with a simple and entirely explicit example %(Example~\ref{exm:a,b-ladder}) 
where the representation by equivalence classes given at the end of \S\ref{sec:Operator-theoretic-interpretation-of-bdG} takes on an especially concrete and visual form. Moreover, the computations can be completed without the direct construction of \Schw, $\Schw'$, or any discussion of $L^2(\Schw',\prob)$; we can obtain the boundary simply by constructing certain functions on the network. We feel this is an especially nice feature of our approach.
%In this section, we introduce the most basic family of examples that illustrate our technical results and exhibit the properties (and support the types of functions) that we have discussed above. %After presenting some basic examples, we prove some theorems regarding the properties of these examples. 

\begin{exm}[One-sided infinite ladder network]\label{exm:a,b-ladder}
  Consider two copies of the nearest-neighbour graph on the nonnegative integers $\bZ^+$, one with vertices labelled by $\{x_n\}$, and the other with vertices labelled by $\{y_n\}$. Fix two positive numbers $\ga > 1 > \gb > 0$. In addition to the edges $\cond_{x_n, x_{n-1}} = \ga^n$ and $\cond_{y_n, y_{n-1}} = \ga^n$, we also add ``rungs'' to the ladder by defining $\cond_{x_n,y_n} = \gb^n$:
  \linenopax
  \begin{equation}\label{eqn:exm:one-sided-ladder-model}
    \xymatrix{
    *+[l]{x_0} \ar@{-}[r]^{\ga} \ar@{-}[d]_{1}
      &  x_1 \ar@{-}[r]^{\ga^2} \ar@{-}[d]_{\gb} 
      &  x_2 \ar@{-}[r]^{\ga^3} \ar@{-}[d]_{\gb^2} 
      &  x_3 \ar@{-}[r]^{\ga^4} \ar@{-}[d]_{\gb^3}
      &\dots \ar@{-}[r]^{\ga^n}
      &  x_n \ar@{-}[r]^{\ga^{n+1}} \ar@{-}[d]_{\gb^n} & \dots \\
    *+[l]{y_0}\ar@{-}[r]^{\ga} 
      &  y_1 \ar@{-}[r]^{\ga^2} 
      &  y_2 \ar@{-}[r]^{\ga^3} 
      &  y_3 \ar@{-}[r]^{\ga^4} 
      &\dots \ar@{-}[r]^{\ga^n}
      &  y_n \ar@{-}[r]^{\ga^{n+1}} 
      & \dots
    }
  \end{equation}
  
  This network was suggested to us by Agelos Georgakopoulos as an example of a one-ended network with nontrivial \Harm. The function $u$ constructed below is the first example of an explicitly computed nonconstant harmonic function of finite energy on a graph with one end (existence of such a phenomenon was proved in \cite{CaW92}). Numerical experiments indicate that this function is also bounded (and even that the sequences $\{u(x_n)\}_{n=0}^\iy$ and $\{u(y_n)\}_{n=0}^\iy$ actually converge very quickly), but we have not yet been able to prove this. Numerical evidence also suggests that \Lap is not essentially self-adjoint on this network, but we have not yet proved this, either.
  
  This graph clearly has one end. We will show that such a network has nontrivial resistance boundary if and only if $\ga > 1$ and in this case, the boundary consists of one point for $\gb=1$, and two points for \gb such that $(1 + \frac1\ga)^2 < \ga / \gb^2$. It will be made clear that the paths $\cpath_x=(x_1,x_2,x_3,\dots)$ and $\cpath_y=(y_1,y_2,y_3,\dots)$ are equivalent in the sense of Definition~\ref{def:harmonic-equivalence-of-paths} if and only if $\gb = 1$.
  
  For presenting the construction of $u$, choose $\gb < 1$ satisfying $4\gb^2 < \ga$ (at the end of the construction, we explain how to adapt the proof for the less restrictive condition $(1 + \frac1\ga)^2 < \ga/\gb^2$). We now construct a nonconstant $u \in \Harm$ with $u(x_0) = 0$ and $u(y_0)=-1$. If we consider the flow induced by $u$, the amount of current flowing through one edge determines $u$ completely (up to a constant). Once it is clear that there are two boundary points %$\gb_0$ and $\gb_1$ 
  in this case, it is clear that specifying the value of $u$ at one %$u(\gb_1)$ 
  (and grounding the other) %$u(\gb_0)=0$) 
  determines $u$ completely.
  
  Due to the symmetry of the graph, we may abuse notation and write $n$ for $x_n$ or $y_n$, and $\check n$ for the vertex ``across the rung'' from $n$. For a function $u$ on the ladder, denote the horizontal increments and the vertical increments by  
  \linenopax
  \begin{align*}%\label{eqn:}
    \gd u(n) := u(n+1) - u(n)
    \qq\text{and}\qq
    \gs u(n) := u(n) - u(\check n),
  \end{align*}
  respectively. 
  Thus, for $n \geq 1$, we can express the equation $\Lap u(n) = 0$ by
  \linenopax
  \begin{align*}%\label{eqn:}
    \Lap u(n) = \ga^n \gd u(n-1) - \ga^{n+1} \gd u(n) + \gb^n \gs u(n) = 0,
  \end{align*}
  which is equivalent to 
  \linenopax
  \begin{align*}%\label{eqn:}
    \gd u(n) = \frac1\ga \gd u(n-1) + \frac{\gb^n}{\ga^{n+1}} \gs u(n).  
  \end{align*}
  Since symmetry allows one to assume that $u(\check n) = 1- u(n)$, we may replace $\gs u(n)$ by $2u(n)+1$ and obtain that any $u$ satisfying
  \linenopax
  \begin{align}\label{eqn:u-harmonic-on-a,b-ladder}
    u(n+1) = u(n) + \tfrac{u(n)-u(n-1)}\ga + \tfrac2\ga\left(\tfrac\gb\ga\right)^n u(n)
     + \tfrac1\ga \left(\tfrac\gb\ga\right)^n
    %\frac1\ga \gd u(n-1) + \frac{\gb^n}{\ga^{n+1}}(2u(n)-1) + u(n)
  \end{align}
  is harmonic. It remains to see that $u$ has finite energy.
  
  Our estimate for $\energy(u) < \iy$ requires the assumption that $\ga > 4\gb^2$, but  numerical computations indicate that $u$ defined by \eqref{eqn:u-harmonic-on-a,b-ladder} will be both bounded and of finite energy, for any $\gb < 1 < \ga$. First, note that $u(1) = \frac1\ga$ and so an immediate induction using \eqref{eqn:u-harmonic-on-a,b-ladder} shows that $\gd u(n) = u(n+1)-u(n)>0$ for all $n \geq 1$, and so $u$ is strictly increasing. Since $\gb < 1 < \ga$, we may choose $N$ so that 
  \linenopax
  \begin{align*}%\label{eqn:}
    n \geq N \q\implies\q \left(\frac\gb\ga\right)^n < \frac{\ga-1}2.
  \end{align*}
  Then $n \geq N$ implies 
  \linenopax
  \begin{align}\label{eqn:ladder(u)-subexponential}
    u(n+1) \leq 2u(n) + \frac1\ga,
  \end{align}
  by using \eqref{eqn:u-harmonic-on-a,b-ladder} and the fact that $u(n)$ is increasing and $\frac\gb\ga < 1$.
  Now use \eqref{eqn:u-harmonic-on-a,b-ladder} to write
  \linenopax
  \begin{align*}%\label{eqn:ladder(du)-iterated-out}
    \gd u(n)
    &= \tfrac1\ga (\gd u)(n-1) + \left(\tfrac2\ga u(n) + \tfrac1\ga\right) \left(\tfrac\gb\ga\right)^n \notag \\
    &= \tfrac1{\ga^n} (\gd u)(0) + \sum_{k=0}^{n-1} \tfrac1{\ga^k} \left(\tfrac2\ga u(n-k) + \tfrac1\ga\right) \left(\tfrac\gb\ga\right)^{n-k} \notag \\
    &= \tfrac1{\ga^{n+1}} + \tfrac{\gb(1-\gb^n)}{\ga^{n+1}(1-\gb)} 
     + \tfrac2{\ga^{n+1}} \sum_{k=1}^{n} \gb^k u(k),
  \end{align*}
  where the second line comes by iterating the first, and the third by algebraic simplification. Applying the estimate \eqref{eqn:ladder(u)-subexponential} gives 
  \linenopax
  \begin{align*}%\label{eqn:}
    2 \sum_{k=1}^{n} \gb^k u(k)
    &\leq 2^2 \sum_{k=1}^{n} \gb^k u(k-1) + \tfrac2\ga \sum_{k=1}^{n} \gb^k
        = 2^2 \sum_{k=2}^{n} \gb^k u(k-1) + 2\tfrac\gb\ga \cdot \tfrac{1-\gb^n}{1-\gb},
  \end{align*}
  and iterating gives
  \linenopax
  \begin{align}\label{eqn:ladder(du)-final-est}
    \gd u(n)
    \leq \frac1{\ga^{n+1}} \left(1 + \frac{\gb(1-\gb^n)}{1-\gb} + \frac{(2\gb)^n}\ga + 2\frac\gb\ga \sum_{k=0}^{n-1} 2^k \frac{\gb^k-\gb^n}{1-\gb}\right).
  \end{align}
  Now the energy $\energy(u) = \sum_{n=0}^\iy \ga^{n+1} \left(\gd u(n)\right)^2$ can be estimated by using \eqref{eqn:ladder(du)-final-est} as follows:
  \linenopax
  \begin{align*}%\label{eqn:}
    \energy(u) 
    &\leq \sum_{n=0}^\iy \frac1{\ga^{n+1}} \left(1 + \frac{\gb(1-\gb^n)}{1-\gb} + \frac{(2\gb)^n}\ga + \frac{2\gb + 2\gb^{n+1} - 2^{n+2}\gb^{n+1} - 2^2\gb^{n+2} + (2\gb)^{n+2}}{\ga (1-\gb) (2\gb-1)}\right)^2
  \end{align*}
  and the condition $\ga > 4\gb^2$ ensures convergence.
  
  Note that this computations above can be slightly refined: instead of $\ga > 4\gb^2$, one need only assume that $\ga > (1 + \frac1\ga)^2\gb^2$. Then, fix $\ge>0$ for which $\ga / \gb^2 > (1 + \frac1\ga)^2 + \ge$ and choose $N$ so that $n \geq N$ implies $\left(\gb/\ga\right)^n < 1 + \frac1\ga + \ge(1 + 2\ga + \ga\ge)$. Then the calculations can be repeated, with most occurrences of $2$ replaced by $1+\frac1\ga+\ge$. 
\end{exm}

\begin{remark}\label{rem:comp-ladder-example-to-Z1}[Comparison of Example~\ref{exm:a,b-ladder} to the 1-dimensional integer lattice]
  In \cite[Ex.~6.3]{DGG}, we showed that the ``nonnegative geometric integers'' network
  \linenopax
  \begin{align*}
    \xymatrix{
      0 \ar@{-}[r]^{\ga} 
      & 1 \ar@{-}[r]^{\ga^2} 
      & 2 \ar@{-}[r]^{\ga^3} 
      & 3 \ar@{-}[r]^{\ga^4} 
      & \dots
    }
  \end{align*} 
  supports a monopole but not a harmonic function of finite energy, for $\ga>1$. These conductances correspond to the biased random walk where, at each vertex, the walker has transition probabilities 
  \linenopax
  \begin{align*}%\label{eqn:}
    p(n,m) = 
    \begin{cases}
      \frac1{1+\ga}, &m=n-1,  \\
      \frac{\ga}{1+\ga}, &m=n+1. 
    \end{cases}
  \end{align*}
  In particular, this is a spatially homogeneous distribution.
  In contrast, the random walk corresponding to Example~\ref{exm:a,b-ladder} has transition probabilities 
  \linenopax
  \begin{align*}%\label{eqn:}
    p(n,m) = 
    \begin{cases}
      \frac1{1 + \ga + \left(\frac\gb\ga\right)^n}, &m=n-1,\\
      \frac{\ga}{1 + \ga + \left(\frac\gb\ga\right)^n}, &m=n+1, \vstr[2.5] \\
      \frac{\left(\gb/\ga\right)^n}{1 + \ga + \left(\frac\gb\ga\right)^n}, &m=\check n. \vstr[2.5]
    \end{cases}
  \end{align*}
  Thus, Example~\ref{exm:a,b-ladder} is asymptotic to the nonnegative geometric integers. 
  
  \pgap
  
  One can even think of Example~\ref{exm:a,b-ladder} as describing the \emph{scattering theory} of the geometric half-integer model, in the sense of \cite{Lax-Phillips}. In this theory, a wave (described by a function) travels towards an obstacle. After the wave collides with the obstacle, the original function is transformed (via the ``scattering operator'') and the resulting wave travels away from the obstacle. The scattering is typically localized in some sense, corresponding to the location of the collision.
  
  To see the analogy with the present scenario, consider the current flow defined by the harmonic function $u$ constructed in Example~\ref{exm:a,b-ladder}, i.e., induced by Ohm's law: $\curr(x,y) = \cond_{xy}(u(x) - u(y))$. With $\act_{_{|\curr|}}(x) := \frac12\sum_{\{z \suth \curr(x,z) > 0\}} |\curr(x,z)|$, this current defines a Markov process with transition probabilities
  \linenopax
  \begin{align*}%\label{eqn:}
    P(x,y) = \frac{\curr(x,y)}{\act_{_{|\curr|}}(x)},
    \q\text{if}\q \curr(x,y)>0,
  \end{align*}
  and $P(x,y) = 0$ otherwise; see \cite{OTERN,RANR}. This describes a random walk where a walker started on the bottom edge of the ladder will tend to step leftwards, but with a geometrically increasing probability of stepping to the upper edge, and then walking rightwards off towards infinity. The walker corresponds to the wave, which is scattered as it approaches the geometrically localized obstacle at the origin.
\end{remark}

\subsection*{Acknowledgements}

We are grateful for helpful (and enjoyable) conversations with Daniel Alpay, Ilwoo Cho, Dorin Dutkay, Agelos Georgakopoulos, Michael Hinz, Keri Kornelson, Paul Muhly, Karen Shuman, Myung-Sin Song and Wolfgang Woess.

\bibliographystyle{alpha}
\bibliography{bdG}

\end{document}